\documentclass[11pt]{amsart}
\usepackage[margin=1.2in]{geometry}

\usepackage{hyperref}

\usepackage[parfill]{parskip}
\usepackage[all]{xy}
\usepackage{latexsym, amssymb, bbm, mathrsfs}
\usepackage{amsthm}
\usepackage{xypic}
\usepackage{xcolor}
\usepackage{graphicx}
\usepackage{amsmath}
\usepackage{enumerate}
\usepackage{arydshln}

\usepackage{comment}

\newtheorem{lem}{Lemma}[section]
\newtheorem{thm}[lem]{Theorem}
\newtheorem{prop}[lem]{Proposition}
\newtheorem{cor}[lem]{Corollary}

\newtheorem{remark}[lem]{Remark}
\newtheorem{assumption}[lem]{Assumption}

\theoremstyle{definition}
\newtheorem{defn}[lem]{Definition}
\newtheorem{example}[lem]{Example}

\newtheorem{question}[lem]{Question}

\newtheorem{exampledefinition}[lem]{Example-Definition}

\newcommand{\ip}[1]{\langle #1 \rangle}
\newcommand{\wt}[1]{\widetilde{#1}}
\newcommand{\wh}[1]{\widehat{#1}}

\newcommand{\ov}[1]{\overline{#1}}

\newcommand{\fmod}[2]{#1\!\operatorname{-}\!#2}

\newcommand{\bdf}{\begin{defn}}
\newcommand{\edf}{\end{defn}}

\newcommand{\bthm}{\begin{thm}}
\newcommand{\ethm}{\end{thm}}

\newcommand{\blem}{\begin{lem}}
\newcommand{\elem}{\end{lem}}

\newcommand{\bcor}{\begin{cor}}
\newcommand{\ecor}{\end{cor}}

\newcommand{\bprop}{\begin{prop}}
\newcommand{\eprop}{\end{prop}}

\newcommand{\brmk}{\begin{remark}}
\newcommand{\ermk}{\end{remark}}

\newcommand{\bpf}{\begin{proof}}
\newcommand{\epf}{\end{proof}}

\newcommand{\bex}{\begin{example}}
\newcommand{\eex}{\end{example}}

\newcommand{\beq}{\begin{equation}}
\newcommand{\eeq}{\end{equation}}

\numberwithin{equation}{section}

\def\co{\colon\thinspace}

\def\nin{\noindent}

\def\sL{\mathscr{L}}

\def\sR{\mathscr{R}}
\def\sS{\mathscr{S}}
\def\sT{\mathscr{T}}

\def\sX{\mathscr{X}}
\def\sY{\mathscr{Y}}

\def\C{\mathbb{C}}

\def\Q{\mathbb{Q}}
\def\R{\mathbb{R}}
\def\Z{\mathbb{Z}}
\def\T{\mathbb{T}}

\def\bl{\mathbbm{l}}
\def\bm{\mathbbm{m}}

\def\bfk{\mathbf{k}}

\def\fm{\mathfrak{m}}

\def\cA{\mathcal{A}}
\def\cB{\mathcal{B}}
\def\cC{\mathcal{C}}
\def\cD{\mathcal{D}}
\def\cE{\mathcal{E}}
\def\cF{\mathcal{F}}
\def\cG{\mathcal{G}}
\def\cH{\mathcal{H}}
\def\cI{\mathcal{I}}
\def\cJ{\mathcal{J}}

\def\cL{\mathcal{L}}
\def\cM{\mathcal{M}}

\def\cO{\mathcal{O}}
\def\cP{\mathcal{P}}
\def\cQ{\mathcal{Q}}
\def\cR{\mathcal{R}}
\def\cS{\mathcal{S}}
\def\cT{\mathcal{T}}
\def\cU{\mathcal{U}}

\def\cY{\mathcal{Y}}

\def\CO{\mathcal{C}\mathcal{O}}
\def\OC{\mathcal{O}\mathcal{C}}
\def\fuk{\mathcal{F}uk}
\def\ov{\overline}

\def\a{\alpha}

\def\e{\epsilon}

\def\la{\lambda}

\def\vp{\varphi}

\def\w{\omega}
\def\z{\zeta}
\def\vp{\varphi}

\def\nin{\noindent}

\def\la{\langle}
\def\ra{\rangle}

\def\xkm2{\overline{X}_{k-2}}
\def\HHGGB{HH_*^G(G\cB)}

\def\Dij{\cD_{i,j}}
\def\Dijp{\cD_{i,j}^\pm}
\def\ij{{i,j}}
\def\ol{\overline}

\def\fix{{fix}}
\def\st{{strict}}
\def\strict{{strict}}

\def\Tan{T^{an}}
\def\Gla{\Gamma^{\bar\alpha}_{l}}

\begin{document}
\title{Equivariant split generation and mirror symmetry of special isogenous tori}



\author[Wu]{Weiwei Wu}
\address{Centre de recherches mathematiques, Universite de Montreal}
\email{wuweiwei@crm.umontreal.ca}

\date{\today}


\maketitle

\begin{abstract}
We prove a version of equivariant split generation of Fukaya category when a symplectic manifold admits a free action of a finite group $G$.  Combining this with some generalizations of Seidel's algebraic frameworks from \cite{Seidelbook}, we obtain new cases of homological mirror symmetry for some symplectic tori with non-split symplectic forms, which we call \textit{special isogenous tori}.  This extends the work of Abouzaid-Smith \cite{AS10}.  We also show that derived Fukaya categories are complete invariants of special isogenous tori.
\end{abstract}

\section{Introduction}

In this paper we investigate the relation between the Fukaya
category of a symplectic manifold and its finite coverings.  This aspect of Fukaya categories has been considered in many different perspectives.  For example, in \cite{Se03} Seidel related an equivariant Fukaya category of a branched cover over a quartic surface to the one downstairs, which led to the first proof of homological mirror symmetry of complex dimension higher than $1$.  This idea was exploited further in \cite{SeiG2}\cite{Sher11}\cite{CH13}\cite{CHL13}\cite{Sher14} for many other instances of homological mirror symmetry.
Another closely related result was given in \cite{RS13}, where Alex Ritter and Ivan Smith showed that
for a finite covering $\pi:X\rightarrow \overline X$ with deck
transformation $G$, if $\overline\cB$ is a collection of Lagrangian
branes that split generates the Fukaya category $\fuk(\overline X)$, then their lifts also split generates $\fuk(X)$.


We would like to address in general the other direction of Ritter-Smith's result: let $\pi: X\rightarrow\ol X$ be a finite covering, if $\cB$ split generates $\fuk(X)$, do their images under the projection $\pi$ also split generate $\fuk(\overline X)$?  The answer is both yes and no: one immediately notices that not all images in $\pi(\cB)$ are guaranteed to be embedded thus are not even objects of Fukaya category of $\ol X$ in the common definition.  However, it is conceivable that, if one is willing to include immersed objects in the formulation of the Fukaya category, the result should still hold by establishing an appropriate version of Abouzaid's generation result \cite{AbGen}.

The approach we adopt in this article is technically simpler. We first define a version of equivariant Fukaya category $\fuk(X)^G$.  The particular formulation we use here is very close to that in \cite{CH13}, to which we will compare in Section \ref{s:defFuk}.  This category, as promised, a technical replacement of the immersed Fukaya category of $\ol X$, in which the extra immersed objects are replaced by $G$-orbits of their lifts in $X$, hence simplifying the situation.

The main tool of our study of the equivariant Fukaya categories is the following $G$-equivariant version of Abouzaid's generation criterion.

\bthm[Corollary \ref{c:gen}]\label{t:gen} Let $\cB$ be a subcategory of $\fuk(X)$. If the open-closed string map
$\OC: HH_*(\cB)\rightarrow HF^*(X)$ hits the identity, then $G\cdot
B$ split generates $\fuk(X)^G$.

\ethm

Following \cite{AbGen}, we say the subcategory $\cB$ \textit{resolves the diagonal} if it satisfies the condition in Theorem \ref{t:gen}.
To relate this equivariant version of Fukaya category to the usual one, we consider the following fully faithful functor from $\fuk(\overline X)$ to the $G$-invariant
category $\fuk(X)^G$.

\bthm[Theorem \ref{t:gen2}]\label{t:functor} There is a transfer functor $\cT:\fuk(\overline X)\rightarrow \fuk(X)^G$ which is full and faithful. \ethm

The functor $\cT$ is naturally defined using the classical transfer
map for covering spaces in the Floer context. By definition, it is easy to see that $\cT$ is an equivalence in the following circumstance:

\bcor[Corollary \ref{c:downstairsGen}] If $\cB\in Ob(\fuk(X))$ resolves the diagonal, and $\pi(L)$ is embedded for any $L\in\cB$, then the collection
$\pi(\cB)$ split generates $\fuk(\ov X)$.

\ecor

Our next task is to understand some applications of this equivariant split generation mechanism.  We are first interested in a general algebraic reduction scheme of a mirror functor when the symplectic side is equipped with a free finite $G$-action.  We set up the problem in a rather algebraic manner as follows.

Given an equivalence functor $\cF:\cC\rightarrow\cD$ between triangulated categories, where $\cC$ is endowed with a strict $G$-action for a finite group $G$.  $\cD$ does not inherit a natural strict $G$-action from $\cF$ in general, but only a \textit{coherent} $G$-action.  We showed that, when the $G$-action can be lifted up to $A_\infty$-level in an appropriate sense, then one may reduce $\cF$ to a fully faithful functor $\cF^\fix$ from $\cC^\fix$ to $\cD^\fix$, while the two fixed (non-full) subcategories and $\cF^\fix$ are again shown to be triangulated.  Also, to obtain a more natural set up to study $\cD^\fix$ when $\cD$ is only endowed with a coherent $G$-action, we propose to consider a \textit{strictification of }$G$-action $\cD^\st$ of $\cD$ in Section \ref{s:strictification}.  We proved the existence of a strictification model for any \textit{finitely generated groups} (although this will not be used in the rest of our paper).  These are again reminiscent of an idea of Seidel in \cite[(14b)]{Seidelbook}, where he considers the case when $G=\Z/2$.  Combining the results from the first part, given a mirror functor $\mathfrak{m}:D^\pi\fuk(X)\rightarrow D^b Coh(X^\vee)$ with a finite free $G$-action on $X$, we have fully faithful functors $\mathfrak{m}^\fix:D^\pi\fuk(X/G)\rightarrow D^b(X^\vee)^\fix$ and $\bar{\mathfrak{m}}^\fix:D^\pi\fuk(X/G)\rightarrow (D^b(X^\vee)^\st)^\fix$.

We then turn to a new case of homological mirror symmetry following a suggestion of Paul Seidel.  We call a symplectic form $\w_{lin}$ of $\R^{2n}$ \textit{linear}, if its coefficients are constant everywhere.  For any linear symplectic forms and a full lattice $\Gamma<\R^{2n}$, $(\R^{2n},\w_{lin})/\Gamma$ is a smooth $\T^{2n}$ endowed with a quotient symplectic form, which by abuse of notation will be denoted as $\w_{lin}$ again.  Such symplectic forms on $\T^{2n}$ will also be called \textit{linear}.

Note that the symplectomorphism type of linear symplectic forms on $\T^{2n}$ are determined completely by linear algebra.  Namely, a symplectomorphism between two linear symplectic forms induces a linear map on $H^1(\T^{2n},\R)$.  Fixing an integral basis in the first cohomology groups, such a map belongs to $GL(2n,\Z)$, which shows that the cohomology classes of the two linear symplectic forms are congruent (as an anti-symmetric bilinear form) by such a matrix.  Since linear symplectic forms are completely determined by the cohomology classes, one may indeed choose a symplectomorphism which lifts to a linear transformation on the universal cover.

The main computation concerns the mirror symmetry of a special type of linear symplectic tori, which we call the \textit{special isogenous symplectic tori}, denoted $T(\bar\alpha)_{\bar l}$.  Roughly speaking, these are tori finitely covered by split symplectic tori in a specific way.  On the $B$-side, we consider abelian varieties $A(\bar\alpha)_{\bar l}$ over the Novikov field, which are again isogenous in a specific way to split analytic tori.  We refer the reader to Definition \ref{d:symptori} and Example \ref{ex:si} for the precise definitions.  Our main homological mirror symmetry statement reads:

\bthm\label{t:mirror} $D^\pi(\fuk(T(\bar\alpha)_{\bar l}))$ is equivalent to $D^b(A(\bar\alpha)_{\bar l})$.
\ethm

The reconstruction theorem due to Bondal and Orlov \cite{BO01} shows that the derived category of an algebraic variety with ample (anti-)canonical line bundle completely determines the variety.  In other words, the derived category is a complete invariant of varieties of this sort.  In contrast, this is not the case for abelian varieties \cite{Mu81}.  Moreover, Polishchuk \cite{Po96} and Orlov \cite{Or02} gave an explicit criterion for two abelian varieties to be derived equivalent over a field $char(\mathbbm{k})=0$.   As the mirror of derived categories in algebraic geometry, the reconstruction theorem of a Fukaya category still seems to be out of reach currently, however, it is still curious how far the Fukaya category is from a complete invariant.  With the mirror symmetry theorem \ref{t:mirror}, one verifies the derived Fukaya category is a complete invariant of special isogenous tori.


\bthm\label{t:distinguish} Two special isogenous tori are symplectomorphic if and only if they have equivalent derived Fukaya categories.
\ethm

The proof crucially relies on Orlov's result, however, the verification of Orlov's condition is far from straightforward: we will need to involve flavors of rigid analytic geometry, which we will recall in Section \ref{s:analSI} and \ref{s:dual}.  We further propose the following question:

\begin{question} Is the derived Fukaya category a complete invariant of linear symplectic tori?
\end{question}

 An affirmative answer to this question should be useful for distinguishing symplectic manifolds of shapes $T(\bar\alpha)_{\bar l}\times M$, where the elementary linear algebra method would no longer work.  This will be the topic of a forthcoming work.

\textbf{Notation:} Throughout $G$ will be a finite group acting
freely on a symplectic manifold $M$ unless otherwise specified. When $L\subset M$ is Lagrangian submanifold,

\begin{itemize}

\item $\ol M=M/G$, and $\ol L=L/G$ if $G$ preserves the Lagrangian submanifold $L\subset M$;

\item $G_L=\{g\in G: g(L)=L\}$ is the isotropy group of $L$;

\item $GL=G\cdot L=\bigcup_{g\in G} g\cdot L;$

\item when $x\in CF^*(L_0,L_1)$, we denote
$$G\cdot x=\bigoplus_{g\in G} g\cdot x\in \bigoplus_{g\in G}CF^*(gL_0,gL_1);$$

\item given a strict/coherent $G$-action on a category $\cC$, we denote $\cC^\fix$ or $\cC^G$ as the subcategory consisting of invariant objects and morphisms. The same applies to a naive $G$-action  on an $A_\infty$-category.

\item the universal Novikov ring is
\beq\label{e:novikov} \Lambda_R=\{\sum_{i=0}^\infty a_iT^{\lambda_i},a_i\in R, \lambda_i\rightarrow +\infty,\lambda_i\in\R\}
\eeq
for any commutative ring $R$.
\end{itemize}

\textbf{Standing Assumption:} To simplify the technicality of the paper, we assume throughout that

\begin{enumerate}[(1)]
\item All Lagrangian submanifolds we consider are spin.  We also require that they either bounds no holomorphic disks for generic choice of compatible almost complex structures; or monotone, i.e.

\beq\label{e:MonoCond}[\w]=\beta\cdot[c_1(M,L)], \hskip3mm\beta>0\eeq

and

\beq\label{e:spin} w_1(L)=0=w_2(L).
\eeq

See more discussions on the monotonicity condition from 4.1.2.

\item $gcd(ord(G), char(R))=1$.  When a $\Z/N$-grading is considered for a symplectic/Lagrangian manifold, $gcd(ord(G), N)=1$ (see Section \ref{s:defFuk} for definition and discussions on gradings).
\end{enumerate}

\noindent\textbf{Acknowledgements:}  The author is particularly grateful to Mohammed Abouzaid for very informative discussions on \cite{AS10} and many other aspects of homological mirror symmetry, which was crucial for the author to initiate this project during the ``Workshop on Moduli Spaces of Pseudo-holomorphic Curves I" in Simons Center; and to Paul Seidel, who suggested considering the special isogenous tori as an application of our reduction method.  Discussions with Cheol-Hyun Cho, Octav Cornea, Luis Haug, Richard Hind, Heather Lee, Yanki Lekili, Tian-Jun Li, Cheuk Yu Mak, Dusa McDuff, Egor Shelukhin, Jingyu Zhao have greatly influenced this paper.  Part of this work was completed during the author's stay in Michigan State University and supported by Selman Akbulut under NSF Focused Research Grants DMS-0244663; Octav Cornea and Francois Lalonde have generously supported trips related to this work; Lingyan Xiao has helped typesetting part of the first draft of this work.  My cordial thanks are due to all of them.

\section{Algebraic Preliminaries}\label{s:algreduction}


The purpose of this section is two-fold: first we would like to recall basic notions of $A_\infty$-categories and fix notations for the rest of the paper.  Then we discuss some purely algebraic results relevant to reducing a mirror functor $\fm:D^\pi(\fuk(M))\rightarrow \cD$ by a finite group action. Here $\cD$ can be any triangulated category, in action it is usually the derived category of coherent sheaves or matrix factorizations of the mirror variety/singularity.

 The reader will note that we almost always focus on the cohomological level hence will mostly only deal with ordinary (triangulated) categories. This is mostly due to the attempt of making our discussion as succinct as possible.  In fact, once the cohomological level is clear, there is a machinery of obstruction theory developed by Paul Seidel \cite[Lecture 14]{Seidelbook2} of upgrading equivariant objects from cohomological level to chain level (weakly equivariant to coherent equivariant in Seidel's terminology), therefore, we will save ourselves from replicating his work.

\subsection{Reminder on $A_\infty$-category}

We collect necessary notion of $A_\infty$-categories we will need, mostly from \cite{Seidelbook} without proofs. Interested readers are
referred thereof for a systematic treatment on the topic.

\bdf\label{d:Ainfty} Fix an arbitrary field $\bfk$. An $A_\infty$-category consists of the following data:

\begin{enumerate}
\item a set of objects $Ob(\cA)$,
\item a graded $\bfk$-vector space $hom_{\cA}(X_0, X_1)$ for each $X_0, X_1 \in$ $Ob(\cA)$ ,
\item a $\bfk$-linear composition maps, for each $d \geq 1$,
\[
\mu_{\cA}^d: hom_{\cA}(X_{d-1},X_d) \otimes \cdots \otimes hom_{\cA}(X_0,X_1) \longrightarrow hom_{\cA}(X_0,X_d)[2-d],
\]
which satisfies quadratic equations
\[
\sum_{m,n} (-1)^{\maltese_n} \mu_{\cA}^{d-m+1} (a_d,\ldots,a_{n+m+1}, \mu^m_{\cA}(a_{n+m},\ldots,a_{n+1}),a_n \ldots, a_1) = 0
\]
with $\maltese_n = \sum_{j=1}^n |a_j|-n$ and where the sum runs over all possible compositions: $1 \leq m \leq d$, $0\leq n \leq d-m$.
\end{enumerate}
\edf

In particular, $hom_{\cA}(X_0,X_1)$ is a cochain complex with differential $\mu^1_{\cA}$; the cohomological category $H(\cA)$ has the same objects as $\cA$ but morphism groups are the cohomologies of these cochain complexes.  In this case, the natural composition maps inherited from $\mu^2$ are associative.

\bdf\label{d:functor} An $A_{\infty}$-functor $\cF: \cA \rightarrow \cB$ between $A_{\infty}$-categories $\cA$ and $\cB$ comprises

 \begin{enumerate}

 \item a map $\cF$: $Ob(\cA) \rightarrow$ $Ob(\cB)$,
 \item a sequence of multilinear maps for $d\geq 1$
\[
\cF^d: hom_{\cA}(X_{d-1}, X_d) \otimes \cdots \otimes hom_{\cA}(X_0,X_1) \rightarrow hom_{\cB}(\cF X_0, \cF X_1)[1-d]
\]
satisfying the polynomial equations
\begin{multline*}
\sum_{r} \sum_{s_1+\cdots+s_r=d} \mu_{\cB}^r(\cF^{s_r}(a_d,\ldots,a_{d-s_r+1}),\ldots, \cF^{s_1}(a_{s_1},\ldots, a_{1})) \\
\quad = \ \sum_{m,n} (-1)^{\maltese_n} \cF^{d-m+1}(a_d,\ldots, a_{n+m+1}, \mu_{\cA}^m (a_{n+m},\ldots, a_{n+1}), a_{n},\ldots a_1).
\end{multline*}

\end{enumerate}
\edf

Such an $A_\infty$-functor defines an ordinary functor $H(\cF):H(\cA) \rightarrow H(\cB)$ which takes $[a] \mapsto [\cF^1(a)]$.

\bdf\label{d:QI} If $H(\cF)$ is an isomorphism (resp. full or faithful), we say $\cF$ is a quasi-isomorphism (resp. cohomologically full or faithful).
\edf


For any $A_{\infty}$-category $\mathcal{A}$, one may consider the $A_\infty$-modules over $\cA$.

\bdf\label{d:module} A \textit{left $\mathcal{A}$-module} $\mathcal{M}$ associates to each object $X\in Ob\mathcal{A}$ a graded vector space $\mathcal{M}(X)$, together with maps
\begin{equation} \mu^{k}: hom(X_{k-1}, X_{k})\otimes \ldots\otimes hom(X_0,X_1) \otimes \mathcal{M}(X_0) \to \mathcal{M}(X_k)\end{equation}

which satisfies $A_\infty$ relations from Definition \ref{d:functor} so that $\cM$ becomes an $A_\infty$-functor from $\cA$ to $Ch$, the dg-category formed by chain complexes over $k$.

\edf

The left $\cA$-modules form a dg-category, which will be denoted $\cA$-mod. Similarly one defines the category of right $\cA$-modules, denoted as mod-$\cA$. There is a $A_\infty$ analogue of Yoneda embedding:

\bdf\label{d:yoneda} Given an object $K\in\mathcal{A}$, we define
 its \textit{Yoneda embedding}, a left module $\mathcal{Y}^l_K$, with
\begin{eqnarray}
\mathcal{Y}^l_K(L) &:= & hom_{\cA}(K,L) \\
\mu^{k}_{\mathcal{Y}^l_K(L)} &:= & \mu^{k+1}_{\cA}.
\end{eqnarray}

In this way Yoneda embedding extends to a cohomologically fully faithful $A_{\infty}$-embedding $\mathcal{Y}^l: \mathcal{A} \to \fmod{\mathcal{A}}{\mbox{mod}}$
and the same holds true for right mod-$\cA$. (for explicit formulae of natural transformations between modules see \cite[Section 2g]{Seidelbook}).
\edf

One of the basic merits $\cA$-module enjoys is the natural triangulated structure. More concretely, if $c\in hom^0_{\cA}(Y_0,Y_1)$ is a cocycle, $Cone(c)$ is an $A_{\infty}$-module defined by
\begin{equation}\label{eqn:cone}
Cone(c)(X) = hom_{\cA}(X,Y_0)[1] \oplus hom_{\cA}(X,Y_1)
\end{equation}
and with operations $\mu^d_{Cone(c)}((b_0,b_1),a_{d-1},\ldots,a_1)$ given by the pair of terms
\[
\left(\mu_{\cA}^d(b_0,a_{d-1},\ldots,a_1), \mu_{\cA}^d(b_1,,a_{d-1},\ldots,a_1)+ \mu_{\cA}^{d+1}(c,b_0,a_{d-1},\ldots,a_1) \right).
\]

In particular one may apply Yoneda embedding and obtain naturally a
\textit{triangulated envelop} of $\cA$ from $\cA$-mod, which is the smallest full subcategory which contains $\cY^l(\cA)$ and are closed under taking cones and applying shift functors. An $A_\infty$-category which is closed under these two operations are called a \textit{triangulated $A_\infty$ category}. A triangulated envelop can also be recast by a concrete construction called the \textit{twisted complex}.

\bdf\label{d:twistedcomplex} A \textit{twisted complex} over $\cA$ is a pair of data $(X,\delta_X)$, so that

\begin{enumerate}[(1)]

\item $X$ is a formal direct sum over finite index set $I$
$$
X = \oplus_{i\in I} V^i \otimes X^i$$
with $\{X^i\} \in$ $Ob(\cA)$ and $V^i$ finite-dimensional graded $\bfk$-vector spaces.

\item $\delta_X$ is a matrix of differentials
$$
\delta_X = (\delta_X^{ji});\qquad \delta_X^{ji} = \sum_k \phi^{jik} \otimes x^{jik}
$$
with $\phi^{jik} \in Hom_k(V^i,V^j)$, $x^{jik} \in hom_{\cA}(X^i,X^j)$ and having total degree $|\phi^{jik}| + |x^{jik}| = 1$. The differential $\delta_X$ should satisfy the two properties
\begin{itemize}
\item $\delta_X$ is strictly lower-triangular with respect to some filtration of $X$;
\item $\sum_{r=1}^{\infty} \mu^r_{\Sigma\cA}(\delta_X,\ldots,\delta_X)=0$.
\end{itemize}
\end{enumerate}

\edf


One observes that twisted complexes form an $A_{\infty}$-category $Tw(\cA)$ (see \cite{Seidelbook}), which is closed under mapping cones and has a natural automorphism by the degree shift functor $\otimes \bfk[1]$. One may show that $Tw(\cA)$ is naturally quasi-isomorphic to the triangulated envelop of the Yoneda modules. For concreteness, we will mostly stick to twisted complexes in this paper. From the triangulated structure of $Tw(\cA)$, one may prove

\blem[\cite{Seidelbook},3.29]\label{l:H0Tw} $H^0(Tw(\cA))$ is an (ordinary) triangulated category.
\elem

Next we will discuss idempotents. Given an additive category $\cC$ and $X\in \cC$, suppose we have an idempotent endomorphism $p\in hom(X,X)$. A triple $(Y,i,r)$ is called an \textit{image of $p$} if the following holds:
$$i\in hom(Y,X),\hskip2mm r\in hom(X,Y),$$
$$ ri=id_Y,\hskip5mm ir=p.$$

$\cC$ is called \textit{idempotent complete} or \textit{Karoubi complete} if every idempotent endomorphism has an image in $\cC$.
An idempotent completion $\Pi\cC$ of $\cC$ can be constructed in a fairly formal way.

\bdf\label{d:idem} The idempotent completion $\Pi\cC$ of an additive category $\cC$ is defined as:

\begin{itemize}
\item $Ob(\Pi\cC)=\{(X,p_X): X\in\cC,\hskip2mm p_X^2=p_X\in hom(X,X)\},$
\item $hom((X_0,p_{X_0}),(X_1,p_{X_1}))=p_{X_1}hom(X_0,X_1)p_{X_0}$
\end{itemize}

\edf

Idempotent completions behave well with respect to the triangulated structure:

\bthm[\cite{BS01}]\label{t:idemtrian} If $\cC$ is triangulated, then $\Pi\cC$ has a natural induced triangulated structure.
\ethm

There is also a notion of \textit{idempotents} in the $A_\infty$ sense.  Let $\cA$ be an $A_\infty$-category as usual.

\bdf\label{d:Aidem} An \textit{idempotent up to homotopy} for an object $Y$ is a non-unital $A_\infty$-functor $\mathscr{P}:\mathbb{K}\rightarrow \cA$ such that $\mathscr{P}(*)=Y$.\edf

 For each idempotent up to homotopy $\mathscr{P}$, one may associate an abstract image of $\mathscr{P}$ which is an $\cA$-module (see \cite[(4b)]{Seidelbook}).  We then say $\cA$ is \textit{split closed} if any such abstract image is quasi-isomorphic to a Yoneda image of an object of $\cA$.

 However, as we mentioned in the beginning of the section, we will only consider idempotents after passing to the homotopy level. This does not lose too much information due to following observation by Seidel.

\blem[\cite{Seidelbook},4c]\label{l:chainidem} Given $A_\infty$-categories $\cA$, $\cB$ and
a cohomologically fully faithful functor $\cF:\cA\rightarrow \cB$. Then $(\cB,\cF)$ forms an $A_\infty$ split-closure of $\cA$ iff $(H^*(\cB), H^*(\cF))$ forms a split-closure of $H^*(\cA)$.

In particular, $H^0(\Pi Tw\cA)$ is equivalent to $\Pi(H^0Tw\cA)$.

\elem

We then define the derived category $D^\pi\cA$ of an $A_\infty$-category $\cA$ as
$$D^\pi\cA:=H^0(\Pi Tw\cA)\cong\Pi(H^0Tw\cA)$$

While the first model is adopted by most of the existing literature, we will also make use of the latter. For situations we will consider, this ambiguity does not incur extra complications.

\subsection{Group actions on a category}\label{s:AlgConcept}

Recall the definition of strict and coherent group action on a category from \cite{Seidelbook}.

\bdf\label{d:actions} Let $G$ be a discrete group, $\cA$ be an additive
category, and $\{T_g\}_{g\in G}$ a set of autoequivalences of $\cA$ such that $T_e=id_{\cA}$.
Then $\{T_g\}_{g\in G}$ forms

\begin{itemize}

\item a \textit{strict} $G$-action if $T_{g_1}\circ T_{g_0}=T_{g_1g_0}$.

\item a \textit{coherent} $G$-action if there is a system of isomorphism of functors $\vp_{g_1,g_0}: T_{g_1}\circ T_{g_0}\xrightarrow{\sim}T_{g_1g_0}$, such that $\vp_{g_1,g_0}$ is the identity when $g_0=e$ or $g_1=e$, and that:

\beq\label{e:coherent}\xymatrix{
T_{g_3}\circ T_{g_2}\circ T_{g_1} \ar[d]_{L_{T_{g_3}}(\vp_{g_1,g_2})} \ar[rr]^{R_{T_1}(\vp_{g_3,g_2})} && T_{g_3g_2}\circ T_{g_1}\ar[d]^{\vp_{g_3g_2,g_1}} \\
T_{g_3}\circ T_{g_2g_1}\ar[rr]^{\vp_{g_3,g_2g_1}} && T_{g_3g_2g_1} }\eeq

\item We also consider the case when $\cA$ is an $A_\infty$ category with a \textit{naive} $G$-action. This means $\{T_g\}_{g\in G}$ is a system of $A_\infty$-autoequivalences so that

\begin{enumerate}[(1)]

\item they form a $G$-action on $Ob(\cA)$,

\item the maps $T_g^1: hom(X,Y)\rightarrow hom(T_g(X),T_g(Y))$ forms a $G$-action,

\item $T_g^i\equiv0$ when $i>1$,


\end{enumerate}



\end{itemize}

\edf


The reason we need both notions of $G$-actions is the following:

\blem\label{l:strictTocoh} Given an equivalence between categories $\cF: \cC\rightarrow \cD$ and a strict $G$-action on $\cC$. Then for any choice of quasi-inverse of $\cG$ of $\cF$, $\{\cF\circ T_g\circ \cG\}_{g\in G}$ forms a coherent $G$-action on $\cD$.\elem

\bpf See \cite[10c]{Seidelbook}. \epf

\bdf\label{d:fix} Let $\cC$ be an additive ($A_\infty$, resp.) category with either a strict or coherent (naive, resp.) $G$-action. The its \textit{fixed part} $\cC^\fix$ denotes the subcategory consisting both objects and morphisms fixed by the $G$-action.
\edf

It goes without saying that the situation of Lemma \ref{l:strictTocoh} models the case of mirror symmetry, where we have constrained our data so that there is a strict action on the $A$-side, and induce a coherent action on the $B$-side.

The proof of Lemma \ref{l:strictTocoh} involves a choice of $\cG$. Much of the theory of $G$-action carries out
without the explicit mention of $\cG$, however, it turns out the invariant part of the $G$-action (the equivariant category) is sensitive to
this particular choice. This means that if we choose the quasi-inverse in an arbitrary way, $\cF$ does not descend to a functor between the invariant part naturally (even when the coherent action is strictified, see Section \ref{s:strictification}). The author do not know if there is a natural way to induce an invariant functor $\cF^\fix$ in such cases. Instead, we impose the following extra condition, which can always be achieved.

\bdf\label{d:admissible} In the situation of Lemma \ref{l:strictTocoh}, a quasi-inverse $\cG:\cD\rightarrow \cC$ is called \textit{admissible} if the following holds:

\begin{enumerate}
\item If $Y\in \cF(\cC^\fix)$, then $\cG(Y)\in\cC^\fix$ and $\cF\circ\cG(Y)=Y$,
\item If $Y_0,Y_1\in\cF(\cC^\fix)$, the pair $\cF: hom(\cG(Y_0),\cG(Y_1))\rightarrow hom(Y_0,Y_1)$ and $\cG:hom(Y_0,Y_1)\rightarrow hom(\cG(Y_0),\cG(Y_1))$ are inverses of each other.
\end{enumerate}

\edf

\blem\label{l:fixtofix} In the situation of Lemma \ref{l:strictTocoh}, given an admissible quasi-inverse $\cG$ of $\cF$, $\cF$ induces a functor $\cF|_{\cC^\fix}:\cC^\fix\rightarrow \cD^\fix$, which is fully faithful.

\elem

\bpf Let $\{T_g\}_{g\in G}$ denote the functors in $G$-action as in Definition \ref{d:actions}. For any $X\in\cC^\fix$, since $\cG\cF(X)\in \cC^\fix$, we have $\cF\circ T_g\circ\cG(\cF X)=\cF \circ\cG(\cF X)=\cF(X)$. The left hand side is precisely the induced coherent action on $\cF(X)$. This proves $\cF(X)\in \cD^\fix$. The property of being fully faithful follows directly from the definition of admissibility (2).

\epf

An important feature of mirror functors is their triangulated structures, hence this property of $\cF^\fix$ is our next topic of study.  To begin with, we assume throughout the rest of this section that:
\begin{itemize}
\item In the situation of Lemma \ref{l:strictTocoh}, $\cC$,$\cD$ and $\cF$ are all triangulated, and the quasi-inverse $\cG$ of $\cF$ is always admissible  (The reader should be reminded that $\cG$ is also triangulated in this case, see \cite[pp. 4]{Happ}).
\item When $G$ is a naive action on an $A_\infty$-category, the $G$-action sends cohomological units to cohomological units (the $G$-action is $c$-unital).
    \end{itemize}

We say that $\cC^\fix$ \textit{inherits a triangulated structure from $\cC$} if it has a triangulated structure so that all exact
triangles in $\cC^\fix$ are triangles in $\cC$.

\blem If $\cC^\fix$ inherits a triangulated structure from $\cC$, then $\cD^\fix$ and $\cF^\fix$ inherits a triangulated structure from $\cD$ and $\cF$.

\elem

\bpf Given $Y_0,Y_1\in \cD^\fix$, there is a triangle in $\cC$ for every $f\in hom(Y_0,Y_1)^\fix$
\beq\label{e:conef} \cG Y_0\xrightarrow{\cG f}\cG Y_1\xrightarrow{i} Cone(\cG f)\xrightarrow{j} \cG Y_0[1].
\eeq

Since $\cG Y_i\in \cC^\fix$ and $\cG f\in hom(\cG Y_0,\cG Y_1)^\fix$ (the latter follows from admissibility), we have a choice of $Cone(\cG f)\in\cC^\fix$ from assumption, and $i$, $j$ are also in the corresponding fixed parts of morphism groups . Applying $\cF$ to \eqref{e:conef}, the following is a triangle in $\cD$ from the admissibility:
\beq\label{e:coneGf}
Y_0\xrightarrow{f}Y_1\xrightarrow{\cF i} \cF(Cone(\cG f))\xrightarrow{\cF j} Y_0[1].
\eeq

But we have seen that $\cF|_{\cC^\fix}:\cC^\fix\rightarrow \cD^\fix$ is a well-defined functor, which implies \eqref{e:coneGf} is indeed a triangle in $\cD$ with all objects and morphisms belonging to $\cD^\fix$. This verifies that $\cD^\fix$ has an inherited triangulated structure from $\cD$. Verification of axioms for triangulated structure on $\cD^\fix$ is straightforward and similar to the above proof.

The triangulated structure of $\cF^\fix$ then follows from that of $\cF$.


\epf

The triangulated structure on $\cC^\fix$ comes for free once we brings in a bit more information from the chain level.  In particular, this always holds when $\cC$ is the derived category of an $A_\infty$-category with a naive $G$-action as the following lemma shows.

\blem\label{l:DAtri} For a triangulated $A_\infty$-category $\cA$ with a naive $G$-action, $D^\pi\cA$ is a triangulated category with a strict $G$-action. Moreover, its fixed part $(D^\pi\cA)^\fix$
is also triangulated.
\elem

This is a mild generalization of a result in \cite{El14} where Elagin proved the result for dg-categories.  Of course, the derived category model we have adopted here is much more elementary, in particular, we did not explicitly involve $A_\infty$-modules hidden in the construction of twisted complexes.

\bpf Recall that using Definition \ref{d:idem}, an object in $D^\pi\cA$ can be represented as

$$\ol Z=(Z=(\bigoplus_{i\in I}V_i\otimes X_i,\delta),p_{Z}),$$

where $p_{Z}\in Hom^0(Z,Z)$ is an idempotent.  Using the notation from Definition \ref{d:actions}, the strict action is given by $\{H^0(T_g)\}$.  Hence, an object of $(D^\pi\cA)^\fix$ is pair consisting a twisted complex and an idempotent, which are both fixed by $G$. Given a $G$-invariant morphism between two objects $f: (Z_1,p_1)\rightarrow(Z_2,p_2)$ in $D^\pi\cA$, if $p_i=id_{Z_i}$, then the usual construction of twisted complex shows the mapping cone $Z_0$ is also an invariant twisted complex. For the general case we argue as in the proof of \cite[Lemma 1.13]{BS01}.



Consider the following diagram of triangles:

\beq\label{e:idem}\xymatrix{Z_1\ar[rr]\ar[d]_{p_1} && Z_2\ar[rr]\ar[d]_{p_2} && Z_0\ar[rr]\ar@{.>}[d]_{t} &&Z_1[1]\ar[d]_{p_1[1]} \\
Z_1\ar[rr] && Z_2\ar[rr] && Z_0\ar[rr] &&Z_1[1]
}\eeq

From the triangulation structure of $H^0Tw\cA$, there is a $t$ which makes the diagram commute. By averaging one may also assume $t$ to be $G$-invariant. Then \cite[Lemma 1.13]{BS01}
shows
$$p_0=t+(t^2-t)-2t(t^2-t)=3t^2-2t^3$$

is an idempotent, by which one may replace $t$ in diagram \eqref{e:idem}. Hence $(Z_0,p_0)$ is the mapping cone of $f$.

Now (TR2) and (TR3) follows from definition and the above averaging trick. This applies to the octahedron axiom as well, but we provide a little more details here. Given a lower cap of the octahedron, one first completes all objects to a $G$-invariant twisted complex by adding another direct summand. This is always possible since we assumed the cohomological unit and all projections involved to be $G$-invariant. One then completes the octahedron by the canonical construction of mapping cones for twisted complexes, which again consists of $G$-invariant twisted complexes. The averaging trick makes all maps on the upper cap $G$-invariant again. Now the rest of the proof follows literally from that of \cite[1.15]{BS01}, by noting that all maps and cones involved are now $G$-invariant, hence this property carries over the whole proof.
\epf


We conclude the general discussions on categorical $G$-actions by bridging the $G$-fixed derived categories and the derived category of a $G$-equivariant $A_\infty$-category, as will be needed in Section \ref{s:defFuk}.


\blem\label{l:dfixfixd} Let $\cA$ be an $A_\infty$-category with a naive $G$-action. Then there is an isomorphism of categories $\cS: \Pi H^0((Tw\cA)^\fix)\hookrightarrow(\Pi H^0 Tw\cA)^\fix=(D^\pi\cA)^\fix$.

\elem

\bpf
The proof is tautological but we include it to dispel possible doubts. Define the functor $\cS:\Pi H^0((Tw\cA)^\fix)\rightarrow(\Pi H^0 Tw\cA)^\fix$ by the natural inclusion.  Given $(Z_i,p_{Z_i})$ with $Z_i\in (Tw\cA)^\fix$, $i=0,1$, where $p_{Z_i}\in Hom^0(Z_i,Z_i)^\fix$ are idempotents. The full and faithfulness of $\cS$ is equivalent to the claim

$$(p_{Z_1}Hom(Z_0,Z_1)p_{Z_0})^\fix=p_{Z_1}Hom(Z_0,Z_1)^\fix p_{Z_0}.$$

But any morphism on the left has the form of a $G$-orbit $G\cdot(p_{Z_1}\circ f\circ p_{Z_1})=p_{Z_1}\circ(Gf)\circ p_{Z_0}$.  On the other hand, for an object $(Z,p)\in \Pi H^0(Tw\cA)$ to be fixed by the $G$-action, $Z$ and $p$ must be both fixed by definition.

\epf



\subsection{Strictification of a coherent $G$-action}\label{s:strictification}

In this section we will explain how to obtain a canonical strict $G$-action model out of a coherent one. The advantage of this model is that the discussion of induced coherent action on the $B$-side becomes more canonical. However, for reducing the mirror functor we still need the admissibility condition on the quasi-inverse $\cG$. Such a strictification model was investigated by Paul Seidel in \cite[14b]{Seidelbook} for the case of $G=\Z/2$. We extend this result to all finitely generated groups using Cayley graph, which is actually much more than what  we need.  One may easily see from our argument that this even works for a broad generality of infinitely generated $G$, as long as one stays in cases that an appropriate version of axiom of choice can be set up, so that the corresponding generalized Cayley graph has a maximal subtree.

\bdf\label{d:CanStrict} Let $\cC$ be a category with a coherent $G$-action $\{T_g, \phi_{g_0,g_1}\}_{g,g_0,g_1\in G}$. Then the \textit{canonical strictification model} $\cC^\st$ is a category equipped with a strict $G$-action, defined by the following data:

\begin{enumerate}[(1)]
\item {\bf Objects}: $\sX\in Ob(\cC^\st)$ has the form $\sX=((\sX_g)_{g\in G},\cF)$, where $\sX_g\in Ob(\cC)$, $\cF=\{\cF_g\}_{g\in G}$. Here $\cF_g\in\bigotimes_{h\in G}hom(\sX_{gh},g\sX_h)$ are isomorphisms. Moreover, we require the compatibility condition

    \beq\label{e:compatibility}
    \left\{\begin{aligned}&\vp^{-1}_{g_1,g_0}\circ\cF_{g_1g_0}=T_{g_1}(\cF_{g_0})\circ\cF_{g_1}&\\
                   &\cF_e=id^{\otimes G}&\end{aligned}\right.
                   \eeq
    Such a system of isomorphisms will be called a \textit{compatible system of isomorphisms}.
\item {\bf Morphisms}: $hom_{\cC^\st}(\sX,\sY)=hom_{\cC}(\sX_e,\sY_e)$ for $\sX=((\sX_g)_{g\in G},\cF^\sX)$, $\sY=((\sY_g)_{g\in G},\cF^\sY)$.

\item {\bf $G$-action}:
\begin{itemize}
\item On the object level, for $s\in G$, $\sX\cdot s=((\wt \sX_g)_{g\in G}, \wt\cF)$. Here $\wt \sX_g=\sX_{gs}$ and $\wt\cF_g\in\bigotimes_{h\in G}hom_{\cC}(\wt \sX_{gh},g\wt \sX_h)=\bigotimes_{h\in G}hom_{\cC}(\sX_{ghs},g\sX_{hs})$ is a shifted copy of $\cF_g$.
\item On the morphism level, if $\Phi\in hom(\sX, \sY)$, then $\Phi\cdot s=(\cF^{\sY}_s)^{-1}\circ T_s(\Phi)\circ\cF^{\sX}_s$.
\end{itemize}
\end{enumerate}

\edf

\brmk\label{r:strict} It is straightforward to check that the $G$-action defined above is strict. The equivalence of $\cC^\st$ and $\cC$ is equally transparent once we construct $\sX$ for each $X\in Ob\cC$ so that $\sX_e=X$ in Proposition \ref{p:strict}. One should notice that the isomorphism type of an object $\sX$ is completely determined by $X_e$. The inheritance of triangulation structure should also be transparent.

However, one also notes that our model has a partially unsatisfactory feature, that we are forced to trade the \textit{left} coherent $G$-action for a \textit{right} strict $G$-action on $\cC^\st$, unless $G$ is abelian. This of course could be remedied in a naive way: we may apply the process again to turn the right action back to a left one when absolutely necessary.
\ermk

Our main proposition of this section is:

\bprop\label{p:strict} For any $G$ finitely generated, if $\cC$ is equipped with a coherent $G$-action, $\cC^\st$ is equivalent to $\cC$.

\eprop

\textbf{A note on notations:} In the proof we will not specify which component of $\cF_g$ is under investigation when it is clear from the context (for example, when its source and target are specified). We will use $t\cdot$ to replace the notation of group action $T_t$ for $t\in G$. Using this notation, $\vp_{t_0,t_1}$ is simply a functor isomorphism from $t_0\cdot t_1\cdot(-)$ to $(t_0t_1)\cdot(-)$. Hence we will also suppress the subscripts of $\vp$ when the context causes no confusions.

\bpf From discussions in Remark \ref{r:strict}, what we need is to construct an $\sX\in Ob(\cC^\strict)$ for any $X\in Ob\cC$, so that $\sX_e=X$. We will show a stronger statement that for any set of objects $\{\sX_g\}_{g\in G}$ satisfying $\sX_g\cong g\sX_e$, there is a compatible system of isomorphisms between these objects.

The first step is to take a set of generators of $G$, $\{t_1,\cdots,t_k\}$, and we consider the Cayley graph $\Gamma$ for this generating set. Recall that

\begin{itemize}
\item $vert(\Gamma)=G$ as a set,
\item there is an oriented edge $e$ with source $s(e)=v_1$ and target $t(e)=v_2$ iff $v_1v_2^{-1}=t_l$ for some $l$. In this case we endow $e$ an extra label $t_l$.
\end{itemize}

Note that our direction of arrows is opposite from the usual notation in Cayley graph.

Replace vertices marked as $g$ by the object $\sX_g$. Our goal is to assign each edge marked by $t_i$ a component of $\cF_{t_i}$, so that the compositions are compatible with \eqref{e:compatibility}. More concretely, if two consecutive edges marked as $t_i$ and $t_j$ are each assigned an isomorphism components of $\cF_{t_i}$ and $\cF_{t_j}$, then \eqref{e:compatibility} uniquely determines a component of $\cF_{t_it_j}$.  Our proposition is equivalent to the assertion that, there is a choice of $\cF_{t_i}$ for edges of the Cayley graph, so that these compositions depends only on the start and end points.  We demonstrate part of the Cayley graph in Figure 1.

\begin{figure}
$$\xymatrix{
 & {\bullet}\ar@{.>}[rrd] & & & & \\
& & &\sX_{t_3g}\ar[r]^{\cF_{t_3}} &\sX_{g} &\\
{\bullet}\ar[ruu] & & & && &{\bullet}\ar@{.>}[llu] \\ & \sX_{t_2t_3g}\ar[rruu]^{\cF_{t_2}}\ar[ul]^{\cF_{t_j}}\ar@{=>}@{.>}@/^5pc/[uurrr]^(.7){\cF_{t_2t_3}} & & &\sX_{t_2'^{-1}t_1'^{-1}t_1t_2t_3}\ar@{.>}[uu]_(.4){\cF_{t_3'}} & \\
\sX_{t_i^{-1}t_1t_2t_3g} & & &\sX_{t_1'^{-1}t_1t_2t_3} \ar[ur]^{\cF_{t_2'}} & &{\bullet}\ar@{.>}[ul]\ar[ruu] \\
 & \sX_{t_1t_2t_3g}\ar[uu]^{\cF_{t_1}}\ar[ul]^{\cF_{t_i}}
 \ar@/_/[rr]_{\cF_{t_l}}\ar[rru]^{\cF_{t'_1}}
\ar@{=>}@{.>}@/^1.5pc/[uuuurr]_{\cF_{t_1t_2}}
\ar@{=>}@{.>}@/_2.5pc/[uurrr]_(.75){\cF_{t_1't_2'}}
\ar@{=>}@{.>}[uuuurrr]^(.6){\cF_{t_1t_2t_3}}
& &{\bullet}\ar[rru]\ar[r] &{\bullet} & }
$$
\caption[caption]{Cayley graph for an object $\sX$. All single arrows are part of the Cayley graph, while the solid arrows belong to $Tr$, and the rest are dotted. Double dotted arrows are not part of the Cayley graph, but elements of isomorphism systems determined by \eqref{e:compatibility}. The top left arrows in the center depicts \eqref{e:comppath} in the actual Cayley grpah. Extra bullets are added for demonstration purposes.}\label{g:Cayley}
\end{figure}
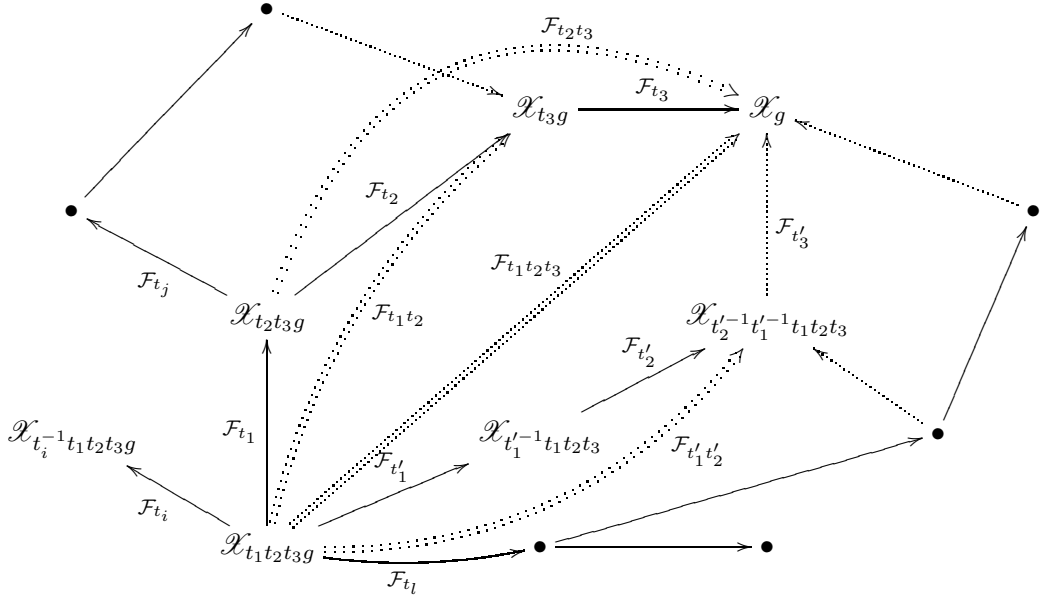

The naive strategy is to start by assigning arbitrary isomorphisms for some edges, then use the above observation to complete the whole compatible system. There are two potential conflicts in completing this process:

\begin{enumerate}[(i)]
\item (compatibility along a path) For triple compositions, $\cF_{t_it_jt_k}$ determined by $\{\cF_{t_it_j}, \cF_{t_k}\}$, and $\{\cF_{t_i}, \cF_{t_jt_k}\}$ must coincide.
\item (compatibility along a loop) When there are more than one oriented path connecting $\sX_{g_0}$ and $\sX_{g_1}$, the isomorphism component of $\cF_{g_0^{-1}g_1}$ determined by the two paths must coincide.
\end{enumerate}

We claim that (i), the compatibility along a path does not impose extra obstructions.

\beq\label{e:comppath}\xymatrix{
\sX_{t_1t_2t_3g}\ar[r]^{\cF_{t_1}} \ar@{.>}[drr] \ar@/_1pc/@{.>}[ddrrr]& t_1\cdot\sX_{t_2t_3g}\ar[r]^{t_1\cF_{t_2}}\ar@{.>}[dd] &t_1\cdot t_2\cdot\sX_{t_3g}\ar[r]^{t_1\cdot t_2\cF_{t_3}}\ar[d]_{\vp} &t_1\cdot t_2\cdot t_3\cdot \sX_g\ar[d]_{\vp}\\
& &(t_1t_2)\cdot\sX_{t_3g}\ar[r]^{(t_1t_2)\cdot\cF_{t_3}}& (t_1t_2)\cdot t_3\cdot\sX_g\ar[d]_{\vp} \\
&t_1\cdot(t_2t_3)\cdot\sX \ar[rr]_{\vp}\ar@{<-}`d[r]`^ul[uurr]_{\vp}[uurr] && (t_1t_2t_3)\cdot\sX_g\\}\eeq

In diagram \eqref{e:comppath}, we have included all possible compositions of relevant isomorphisms, corresponding to three consecutive arrows in the Cayley graph. Arrows in the first row are given arbitrary isomorphism components of $\cF$, and this will determine all the rest of the data in the diagram. To see this yields a commutative diagram, we note first that all the solid arrows form a commutative sub-diagram: there are two loops in question, where the top-right square is commutative since $\vp$ is an isomorphism of functors; the bottom-right diagram is simply the coherence condition of the action. Now by definition, the three dotted arrows are determined by the sub-diagram formed by solid arrows, the commutativity hence follows.


For (ii), we take a connected maximal subtree $Tr$ of $\Gamma$ with the following property:

\begin{itemize}
\item $\sX_e\in Tr$
\item there is an oriented path in $Tr$ starting from $\sX_e$ and ends at $\sX_g$ for all $g\in G$.
\end{itemize}

Note that $Tr$ contains all vertices in $\Gamma$: otherwise, take any missing vertex $\sX_{m}$. $m$ can be decomposed as
$$m=t_{i_1}\cdots t_{i_k}.$$

Let ${\displaystyle l'=\max\{r: 1\le r\le k, \hskip 2mm t_{i_1}\cdot t_{i_2}\dots\cdot t_{i_r}\in Tr\}}$. Adding vertices $\{\Pi_{1\le \nu\le p} t_{i_\nu}\}_{l'\le p\le l}$ and corresponding edges yields a strictly larger subtree than $Tr$.

Now assign arbitrary isomorphism components corresponding to each edge belonging to $Tr$. Then the compatibility on paths proved above yields a system of isomorphisms through \eqref{e:compatibility} on all compositions along $Tr$.

For edges that are not in $Tr$, by definition, adding any of them to $Tr$ yields precisely one loop to the new graph. Hence the prescribed data on $Tr$ and \eqref{e:compatibility} automatically assigns an isomorphism to these edges, and by definition (or repeating the proof of compatibility along paths), this way of assigning isomorphisms yields a compatible system of isomorphisms as desired.

Lastly, notice that we indeed fixed a bit more data than in the proof by requiring $\cF_e\equiv id^{\otimes G}$. But a quick reflection shows this only cause potential problem when there is a loop starting from a vertex and back, which consists a loop that we could deal with as in (ii) in the proof.

\epf

Note that although the embedding $\cC\hookrightarrow \cC^{strict}$ is canonical only up to a choice of compatibility system of isomorphism for each object, it is canonical on $\cC^\fix$ where the system of isomorphisms can be chosen as identity (strictly speaking one needs to require $\vp$ to be identity on $\cC^\fix$, too, but this is only cosmetic). Hence by combining \ref{l:dfixfixd}, \ref{l:fixtofix} from previous sections we obtain the following:

\bprop\label{p:red} If $\cA$ is a triangulated $A_\infty$-category with a naive $G$-action, and there is an triangulated equivalence $\mathfrak{m}:D^\pi\cA\rightarrow \cB$, then the equivalence can be reduced to triangulated fully faithful embeddings

\beq\label{e:reduction}
\begin{array}{l}
\mathfrak{m}^\fix: (D^\pi\cA)^\fix\hookrightarrow \cB^\fix,\\
\bar{\mathfrak{m}}^\fix: (D^\pi\cA)^\fix\hookrightarrow (\cB^{strict})^\fix.
\end{array}
\eeq

\eprop

\textit{Question: When are $\mathfrak{m}^\fix$ and $\bar{\mathfrak{m}}^\fix$ exact equivalences?}

Abstractly, this question is related to the equivariant theory recently developed by Paul Seidel in \cite{SeAnEq}\cite{Seidelbook2}.
Using terminologies from his work, the potential failure for $(D^\pi\cA)^{fix}\rightarrow \cB^{fix}$ to be essentially surjective lies in that whether one could find a weakly equivariant object with image $Y$ for all $Y\in\cB^\fix$. If this is the case, \cite[Lemma 14.10]{Seidelbook2} provides a way to produce the necessary object in $(D^\pi\cA)^{fix}$. The next step for $\cB^{fix}\hookrightarrow (\cB^{strict})^{fix}$ would probably require additional considerations on obstruction theories for coherent actions.

In this perspective, what we have shown is quite preliminary: given a chain level (for example, a naive one) $G$-action on $\cA$ and a quasi-equivalence to $\cA'$, the strictification model provided an approach so that one could discuss weakly equivariant objects on both sides as a start. However, problems about upgrading the $G$-action or equivariant objects on $\cA'$ to the chain level, as well as the comparison of these objects on $\cA$ and $\cA'$, remain mysterious.
However, we will see later that, in some geometric situations, one may show the equivalence from extracting an isomorphic piece from each side, where the induced $G$-action on $\cB$ is still strict.

\section{Equivariant Fukaya category and split generation}

\subsection{Moduli spaces of bordered holomorphic curves}\label{s:moduli}

We start by considering various moduli spaces of pseudo-holomorphic curves
 involved in our discussions.  Since a comprehensive account on these moduli spaces in rather
general contexts has been carefully written down in
\cite{AbGen}\cite{Gan13}, we only recall relevant notions involved
in our later discussions and refer interested readers to their
detailed expositions.



Let $D^2=\{z\in \C: |z|\le 1\}$. Denote the moduli
space

$$\wt \cD_\ij=\{D^2\backslash \Sigma^+\cup\Sigma^-:\Sigma^+\cup \Sigma^-\subset\partial D^2; |\Sigma^+|=i,|\Sigma^-|=j\},$$

\noindent and

$$\wt \cD_\ij^\pm=\{D^2\backslash \Sigma^+\cup\Sigma^-\cup p^\pm:
\Sigma^+\cup \Sigma^-\subset\partial D^2; |\Sigma^+|=i,|\Sigma^-|=j;
p^\pm\in \text{int}(D^2)\},$$

where $D^2$ is always equipped with the standard complex structure.
In both cases, we order $\Sigma^+$ and $\Sigma^-$ counterclockwisely
as $\{z_1^+,\cdots,z_i^+\}$ and $\{z_1^-,\cdots, z_j^-\}$. To be more
pertinent to our applications, we further restrict our
attention to the following components of the above moduli spaces:

\bdf\label{d:Dij} $\Dij\subset\wt \cD_\ij$, $\Dijp\subset\wt
\cD_\ij^\pm$ are components such that the following additional
restrictions are satisfied:

\begin{enumerate}[(1)]

\item $j=0,1$ or $2$;

\item $\Sigma^-$ lies on the same connected component of $\partial D^2\backslash
\Sigma^+$.

\end{enumerate}

\edf

We will denote the Deligne-Mumford compactification of these two types of moduli spaces as $\overline{\cD}_{i,j}$ and $\overline{\cD}^\pm_{i,j}$,
where the lower strata are stratified by stable disks with more than one components.

Next we consider the moduli space for annuli, denote $A_r=\{1\le|z|\le{}r\}$. and the moduli space $\cC^-_d=\{(A_r,z_0,z_1,...z_d)\vert{}z_0=1,z_1=r,|z_k|=r,\forall k\ge{}2\}$.

The Deligne-Mumford compactified moduli space $\overline\cC^-_d$ has similar boundary strata as $\overline{\cD}_\ij$, i.e bubbling up disks at the boundary $|z|=r$. However, two types of new boundary strata of codimension 1 appears where the annulus breaks into two disks $D_1$ and $D_2$, and either
\begin{enumerate}[(1)]
\item $D_1\in\cD_{d,0}^-$ and $D_2\in\cD_{0,1}^+$
\item\label{2} $D_1\in\cD_{d_{1,2},},D_2\in\cD_{d-d_1+2,1}$
\end{enumerate}

In many occasions below, we will deal with $\Dijp$ and $\cC^-_d$ uniformly. In such cases we will simply use $\sR$ and $\ov\sR$ to denote any of these moduli spaces and their compactifications, respectively.

\subsection{Floer data and Consistent choices}

\bdf\label{d:consistent} A \textbf{$G$-invariant Floer datum} on a disk $S\in \sR$ consists of the following choice on each component:
\begin{enumerate}[(1)]
\item \textit{Strip-like ends and cylindrical ends:} a strip-like end near $z_l^{\pm}$ is a choice of $\e_l^{\pm}:Z_{\pm}\rightarrow S$ , where $Z_+$=$[0,+\infty)\times[0,1]$ , $Z_-=(-\infty,0]\times[0,1]$ and
$\lim\limits_{s\to\pm\infty}\e_l^\pm(s,t)=z_l^\pm$ . A cylindrical end near $p^\pm$ is a choice of $\e_0^\pm:W_\pm\rightarrow{}S$ , where $W_+$=$[0,+\infty)\times{}S^1$ , $W_-=(-\infty,0]\times{}S^1$, and
$\lim\limits_{s\to\pm\infty}\e_0^\pm(s,t)=p^\pm$.

\item \textit{Hamiltonian perturbations:} this is a map $H_S:S\to \cH(M)^{G}$ on each surface defining a Hamiltonian flow $X_S$ depending on points over $S$, where $\cH(M)^G$ denotes the space of Hamiltonian functions on $M$ that are $G$-invariant.  Moreover, $H_S\circ \epsilon_i^\pm(s,t)$ is independent of $s$ for all $i\ge0$.

\item\label{one-form} \textit{Basic one-form:} a basic one-form $\a_S$ satisfies $\a_S\vert_{\partial{}S=0}$ and that $X_S\otimes\a_S=X_{H_S}\otimes{}dt$ at all cylindrical ends or strip-like ends.

\item \textit{Almost complex structures:} this is a map $I_S:S\to{}\cJ(M)^{G}$ whose pullback under $\e_l^\pm$ depends only on $t$, where
$\cJ(M)^{G}$ denotes the $G$-invariant compatible almost complex structure on $(M,\w)$.

 \end{enumerate}
\edf

A special case that we did not include is when we have a disk $D$ with one input $z^+$ and one output $z^-$ . We fix a diffeomorphism
\begin{equation}
\begin{array}{llll}

S_{1,1}=D\setminus\{z^+,z^-\}\simeq{} & \R &\times & [0,1] \\
{} & (s &, & t)
\end{array}
\end{equation}and choose $H_S$ and $I_S$ to be invariant under translations in the $s$-variable for similar construction in \ref{d:consistent}.



With the aim of defining Fukaya category, one needs to restrict the choice of Floer data. Before getting that far, we first explain the notion of
\textbf{universal and consistent choice of Lagrangian labels}. Explicitly, for any $S\in\sR$, one assigns a label $\zeta$ to a connected component of $\partial S\backslash\cup\{z_l\}$ which we denote as $\partial_\z{}S$. One assigns a Lagrangian $L_\z$ to $\partial_\z{}S$, which we call a \textit{Lagrangian label}. Seidel showed in \cite{Seidelbook} that, one could extend a given choice of Lagrangian labels from a Riemann surface $S\in\sS$ to a set of Lagrangian labels for the universal family $\ov\sS$ of Riemannian surfaces over $\ov\sR$, so that it is locally constant over
$\sS$ and compatible (in the most obvious sense) with the gluing maps from lower strata to higher strata.



We now may consider extending the choice of Floer datum from $\sR$ to $\ov\sR$.

\bdf\label{d:univcon} A \textit{universal and consistent choice of $G$-invariant Floer data for} $\sR$ is the following: given a set of Lagrangian labels $\{L_\z\}$, we have a choice of $G$-invariant Floer data for each element $S\in \ov\sR$ which varies smoothly over this compactified moduli space. Floer data between different strata satisfies the following compatibility conditions:
\begin{enumerate}[(1)]
\item On boundary strata the Floer data are conformally equivalent to the product of Floer data over irreducible components.
\item In the coordinates given by the parameters gluing at cylindrical or strip-like ends, the Floer data near the boundary agrees with those obtained by gluing of the boundary strata up to infinite order.
\item At any positive (resp.negative) strip-like end, $\epsilon^+(\e^-,\mbox{resp.})$ with adjacent Lagrangian labels $L_1^+$, $L_2^+$ ($L_1^-$, $L_2^-$ resp.), the Floer data is the chosen one for the strip with corresponding Lagrangian labels.

We will denote the space of $G$-invariant Floer data over a universal family $\sR$ as $\cP^G_{\sR}$ 
In most cases when the universal family is clear from the context, we will suppress the subscript and simply use $\cP^G$.

\end{enumerate}

\edf


Again such choices could be obtained in the non-equivariant $\sR=\cD_{n,1}$ case, according to \cite[(9g,9i)]{Seidelbook}.  See \cite{AbGen} for the adaption when $\sR$ is different from $\cD_{n,1}$. The proof completely carries over to our $G$-equivariant case without modifications. We will show in \ref{s:ET} that
within $\cP^G$ transversality could also be obtain in the monotone cases.

\brmk One notices our definitions for various moduli spaces and consistent choices of data are contained in those appeared in \cite{AbGen}.  Our case is only simpler due to the fact that we do not deal with wrapped Floer cohomology, so that one does not need the delicate rescaling trick, hence discussions on weights can be ignored.

This affects Definition \ref{d:consistent}(\ref{one-form}), where cylindrical ends and strip-like ends need not be distinguished; and Definition \ref{d:univcon}(2), where in the wrapped case the gluing cannot be performed in a naive way as we described (see \cite[Section 6.2]{AbGen}).  We do not encounter such subtlety in this paper, which is a lot more convenient. This in turn means our consistent choice can be understood in the original way that Seidel described in \cite{Seidelbook}.  \ermk
\subsection{Moduli spaces of perturbed holomorphic curves}\label{s:MHC}

Using notations from the previous section, for $S\in\sR$, we consider maps
$u:S\to{}M$, which satisfies the following set of equations:
\beq\label{e:FloerEq}
\left\{ \begin{array}{l}
(du-X_{H_s}\otimes\alpha_S)^{(0,1)}=0 \\
\lim\limits_{s\rightarrow\pm \infty}u\circ\e_0^\pm(s,t)=y^\pm(t),\\
\lim\limits_{s\to\pm\infty}u\circ\e_l^\pm(s,t)=x_l^\pm(t), \hskip 2mm l\ge1,\\
u(\partial_\z{}S)\in L_\z
\end{array}
\right.
\eeq
Let us explain the equations term by term.  We use a universal choice of Floer data associated to a given set of Lagrangian labels $\{L_\z\}$ for $\sR$.  $y^\pm$ is a Hamiltonian orbit of $H_S$ restricted to the cylindrical end $\e_0^\pm$. For a given strip-like end $\e^\pm_l$, $l\ge1$, it has two adjacent-components of $\partial_{\z_0^\pm}S$ and $\partial_{\z_1^\pm}S$ numbered counter-clockwisely when it is a positive end (input) and clockwisely when it is a negative end (output). Then its asymptotic limit $x_l^\pm(-)=\lim_{s\rightarrow\pm\infty}u\circ\e^\pm(s,-)$ is a Hamiltonian chord going from $L_{\z_0^\pm}$ to $L_{\z_1^\pm}$.
These moduli spaces all come with relative orientations with respect to orientations of determinant lines of chords $x_l^\pm$ and orbits $y^\pm$ as long as all Lagrangian labels come with a fixed spin structure, as proved in \cite{Seidelbook} and \cite{FOOO_Book}.

\bdf\label{d:moduli} When $S\in\sR$, for $\sR=\cD_\ij^{(\pm)}$ or
$\ol{\cC}^-_{d}$, $\cM_{\sR}$ denotes the moduli space of  \eqref{e:FloerEq}, and $\ol{\cM}_{\ol\sR}$ as its compactification.
\edf

\subsection{Equivariant Transversality}\label{s:ET}

The transversality problem of \eqref{e:FloerEq} is by now standard under the condition \eqref{e:MonoCond}.  See \cite{BC1}, \cite{BC2} for the treatment for defining Fukaya categories specifically in this setting.
For more general cases one needs virtual perturbations techniques from \cite{FOOO_Book}.

In the equivariant setting with finite $G$-actions, this was addressed in \cite{CH13} again using Kuranishi structures.  However, in our setting when the $G$-action is free, the issue of equivariant transversality can be settled by more or less classical methods such as \cite[Lemma 5.13]{KS02} and \cite[Sectoin (9k)]{Seidelbook}. Although such arguments have appeared in the literature, we feel it instructive to re-iterate part of it here just for the sake of being self-contained, as well as to remove potential doubts, but we will not elaborate all the details.  We will focus on a single Riemann surface case, the family version will follow from adapting arguments in \cite[(9k)]{Seidelbook}

\textbf{Upshot:} given $S\subset \Dij^{(\pm)}$ or $\cC_d^-$ decorated with Lagrangian labels $\{L_\zeta\}$, one may choose equivariant perturbation datum $(H_S, J_S)$ from a generic set $\cP_{gen}\in \cP^G$, so that the moduli problem  \eqref{e:FloerEq} modelled over $S$ achieves transversality with $(H_S, J_S)$.

To this end, one first notices that, the virtual dimension count of bubbles is not affected by a finite free group action (unbranched covers of a sphere are trivial covers), hence the issues of bubbling is taken care of in a completely analogous way as in the non-equivariant case. Then the main point of the proof is that, as in the usual genericity argument, a cokernel element of the Cauchy-Riemann operator $D_{u,J}$ for such a perturbed holomorphic curve $u:S\rightarrow M$ satisfying \eqref{e:FloerEq} gives a nonzero section $Z\in L^q(\Lambda_{S}^{0,1}, u^*TM)$ for some $q>0$, which satisfies a $\bar\partial$-type equation and satisfies:

\beq\label{e:EqTr}\int_{S}\w\langle (\delta Y)^{0,1}, Z\rangle dsdt=0\eeq

\noindent for all $Y\in T_J$, the tangent space at $J$ in the space of domain dependent compatible almost complex structures. We
assume $Z(z_0)\neq0$ for $z_0\in S$.

When $S$ is not a strip, one could first construct $\wt{\delta Y}$ which is not necessarily $G$-invariant but concerntrated near $z_0$ over the domain and $u(z_0)$ over the target, so that the left hand side of \eqref{e:EqTr} is non-zero hence the equality fails. To make the infinitesimal variation $G$-invariant, simply use an averaging process on the target symplectic manifold $M$ to replace $\wt{\delta Y}$ by $\delta Y=\sum_{g\in G}g^*\wt{\delta Y}$ (note that $\delta Y$ still concerntrates near $z_0$ over the domain).

When $S=\R\times [0,1]$ is a strip, there is an additional requirement that $\delta Y$ must be $s$-invariant. In this case we resort to \cite[Lemma 5.12]{KS02} and \cite[Theorem 4.3]{FHS}. A key property is:

\blem[\cite{KS02}, Lemma 5.12(J5)]\label{l:J5} Let $v_1,v_2:\R\times [0,1]\rightarrow M$ be two non-constant solutions of $\partial_s u+J_t\partial_t u=0$. Assume that $v_2$ is not a translate of $v_1$ in $s$-direction. Then for any $\rho>0$ the subset
$S_\rho(v_1,v_2)=\{(s,t)\in\R\times (0,1): v_1(s,t)\notin v_2([-\rho, \rho]\times \{t\})\}$ is open and dense in $\R\times(0,1)$.

\elem

From this we consider sets

\begin{align*}
&R_1(u)=\{(s,t)\in \R\times (0,1): du(s,t)\neq0, u(s,t)\notin u(\R\backslash\{s\}, t), u(s,t)\neq x_\pm\},\\
&R_2(u)=\{(s,t)\in \R\times (0,1): u(s,t)\neq g(x_\pm), \forall g\in G, g\neq e\},\\
&R_3(u)=\{(s,t)\in \R\times (0,1): u(s,t)\notin g(u(\R,t)), \forall g\in G, g\neq e\}, \end{align*}

\nin where $x_\pm$ are the limits for $u(s,t)$ when $s\rightarrow\pm\infty$. $R_1$ and $R_2$ are both open and dense (see for example \cite[Lemma 5.12 (J3)]{KS02}), and $R_3(u)$ is a countable intersection of open dense sets by Lemma \ref{l:J5}. Hence $R(u)=R_1(u)\cap R_2(u)\cap R_3(u)$ is a residual set, and it can indeed be shown to be open as in \cite[Theorem 4.3]{FHS}. With this understood, the averaging process on $s$-invariant $\delta Y$ as above again leads to a contradiction in \eqref{e:EqTr}. This concludes the proof of the equivariant transversality.

As a result of what is explained in this section, we will make no mention to the transversality issue in the rest of the paper and use freely the fact that all moduli problems we encounter satisfies transversality by choosing generic universal and consistent Floer data.

\section{$G$-equivariant generation criterion}\label{s:defFuk}


In this section we define a version of $G$-equivariant
Fukaya category $\fuk(M)^G$. We will first review
ingredients involved in the definition of usual Fukaya category,
which are mostly taken from \cite{Seidelbook},
then explain how to incorporate the $G$-action.

The reader will note that our definition is a mild generalization (with simplifications on some technical points) of that in \cite{Seidelbook} for the case of $G=\Z/2$.  We should also compare our version to another very similar version of $G$-equivariant Fukaya category defined in \cite{CH13}.  For readers' convenience, we list the differences (using their terminology) as follows:

\begin{itemize}

\item In \cite{CH13} the authors discusses general cases of different spin profiles which defines different equivariant Fukaya categories, while we always take the trivial spin profile in their terminology.
\item We restrict ourselves to the monotone cases and the universal Novikov field coefficients, and avoid the use of general $G$-Novikov theory.
\item We allowed intersections between $gL$ and $g'L$ for $g\neq g'\in G$.
\item We include discussions on $\Z/N$-gradings for $N\neq\infty$.
\end{itemize}

In particular, we have restricted our setting convenient for the $G$-equivariant generation, and have paid exclusive attention to issues relevant to passing to the quotient instead of making any attempts to a general equivariant theory.  Presumably, one should still be able to construct a more general version of equivariant split generation after incorporating more ingredients such as non-trivial spin profiles and Kuranishi structures from \cite{CH13}, but we will not discuss this point in the current paper.

\subsection{Review on Fukaya category}

\subsubsection{Technical aspects: gradings and spin structures} Prior to actually defining the Fukaya category, we need to recall two
basic technical ingredients necessary for well-defined degrees and signs appearing in the definition. We will only give brief summaries without proofs on notions involved in the $G$-equivariant variations.


\textit{$\bullet$ Gradings:}

$\Z/N$-gradings of an embedded Lagrangian submanifold for $2\le N\le\infty$ is defined in \cite{SeGraded}. For coherence of notation, in below $\Z/\infty$ will be understood as $\Z$. Recall that if we denote $\sL$ as the total space of the Lagrangian Grassmanian bundle over $(M,\w)$, $\sL^N\rightarrow \sL$ is an $N$-fold cover so that its restriction to $\sL_x$ is the standard $N$-fold cover associated to a preferred generator of $H^1(\sL_x,\Z/N)$. The existence of a $\Z/N$-grading is in turn equivalent to either of the following:

\begin{enumerate}[(a)]
\item Existence of a global Maslov class mod $N$, i.e. $C_N\in H^1(\sL,\Z/N)$, such that $C_N|_{\sL_x}$ is the preferred generator in $H^1(\sL_x,\Z/N)$,

\item $2c_1(M,\w)=0$ in $H^2(M,\Z/N)$.
\end{enumerate}

As a result of $(b)$, there is always a $2$-fold Maslov cover. Assuming the existence of an $N$-fold Maslov covering, we may consider the \textit{$\Z/N$-grading} of a Lagrangian submanifold $L\subset (M,\w)$. This is a lift $gr^N: L\to \sL^N|_L$ from the natural section $L\rightarrow\sL|_L$. The existence of such a grading in turn implies $\sL^N\to \sL$ is a trivial $N$-fold cover when restricted to $TL$. Note that when $L$ is orientable (as is the only case we consider), its orientations gives natural $\Z/2$-gradings, hence we always assume $N\ge 2$.  Following the usual convention, we will call the pair $(L,gr^N)$ for $2\le N\le\infty$ a \textit{$\Z/N$-graded Lagrangian}.

Although the discussion above covers the case of $\Z$-gradings, there is a more explicit way of describing it which is worth recalling.
Take a quadratic volume form $\eta^2\in
(\wedge^n(TM;J)^{\otimes 2})^\vee$, one defines
$\alpha_M:\sL\rightarrow S^1$ as
\beq\label{e:squarephase}\alpha_M(\xi_x)=\eta(v_1\wedge\cdots\wedge
v_n)^2/|\eta(v_1\wedge\cdots\wedge v_n)^2|\in S^1 \eeq
for any basis $\{v_i\}$ of a Lagrangian subspace $\xi_x\subset
TM_x$. This offers an explicit representative of a global Maslov class in $H^1(\sL,\Z)$.
  A grading of a Lagrangian submanifold $L$ is therefore a lift
$\alpha^\#:L\to \R$, so that the composition
$L\xrightarrow[\empty]{\alpha^\#}\R\to\R/\Z=S^1$ coincides with the
restriction of $\alpha_M$ to $TL\subset TM$. The existence of a
grading for $L$ is equivalent to the vanishing of the Maslov class
$\mu_L\in H^1(L;\Z)$.


The existence of gradings has the following implication in Floer theory:
given two $\Z/N$-graded Lagrangians $(L_0,\alpha_0^\#)$, $(L_1,\alpha_1^\#)$,
 any intersection point $x\in L_0\cap
L_1$ obtains an absolute degree $deg_{(L_0,\alpha_0^\#),
(L_1,\alpha_1^\#)}(x)\in \Z/N$ from the extra grading structure as explained
in \cite{SeGraded}. We will abbreviate this as $deg(x)$ when no
confusion occurs. We omit the concrete construction here,
but the following situation will be relevant. Given any
symplectomorphism $f:M\to M$, $(f(L_i),\alpha_i^\#\circ f^{-1})$ are
again two graded Lagrangians, $i=0,1$. Then $deg(f(x))=deg(x)$.

$\bullet$ \textit{Spin structures:}


With a chosen spin structure on
each Lagrangian labels, all moduli spaces involved in the Floer cohomology operations we will consider are all equipped with
a preferred orientation, hence one may talk about signs in a coherent way, see \cite{Seidelbook}. Hence, if one is willing to work over a base ring with $char(R)=2$, this assumption is redundant.

\subsubsection{Definition of Fukaya category $\fuk(M)$.}

With the above technical points understood, we may explain the definition of the Fukaya category $\fuk(M)$ following \cite{Seidelbook}.

$\bullet$ \textit{Objects:}

An object in $\fuk(M)$ is a Lagrangian brane, that is, a triple $L^\#=(L,gr^N_L, spin(L))$, where $L\subset (M,\w)$ is an embedded
Lagrangian submanifold, $gr^N_L$ is a $\Z/N$-grading, and $spin(L)$
a chosen spin structure over $L$. For ease of notation, sometimes the additional data $(gr^N_L, spin(L))$ will be suppressed and we will just refer to a Lagrangian brane $L^\#$ using its underlying Lagrangian $L$ when the context is clear.

The above definition already works for the consideration of Fukaya categories consisting of unobstructed Lagrangians.  Although we have already seen the transversality issues can be taken care of in a standard way for Lagrangians satisfying the Standing Assumption in the introduction, for the Fukaya category to be well-defined in the monotone case, the value of \textit{superpotential} of a monotone Lagrangian in the sense of Fukaya-Oh-Ohta-Ono \cite{FOOO_Book} needs to be recalled.

 Explicitly, the value of superpotential is a Gromov-Witten type invariant which counts the algebraic number£¬ weighted by Novikov coefficients, of holomorphic disks of Maslov index $2$ with the boundary passing through a fixed point of $L$, i.e.
 
 $$m_0(L)=\sum_{\begin{subarray}{l}
 \hskip6mm \mu(\beta)=2\\
 \hskip4mm\beta\in H_2(M,L)\end{subarray}} T^{\w(\beta)}\cdot (ev_0)_*[\cM_{\cD_{0,1}}(\beta)]/[L]\in\Lambda.$$

In \cite[Lemma 3.2]{RS13} and \cite{BC1}\cite{BC2} (for ungraded and $\Z_2$-coefficient case), it was shown that Lagrangians with the same value of $m_0$ consist an $A_\infty$-categories as in the unobstructed case.  In the rest of this paper, our treatment can be made uniform for unobstructed cases and monotone cases only by noting this point. Hence, we will make no explicit mention regarding the $m_0$-value of the Fukaya category under considerations from this point on.

$\bullet$ \textit{Morphisms:}

We define the morphism groups as

$$hom^*(L_0, L_1)=CF^*(L_0, L_1)=\bigoplus_{x\in \cC(L_0,L_1; H_{L_0,L_1})}\Lambda_R\cdot\langle x\rangle.$$
Here $\cC(L_0, L_1;H_{L_0, L_1})$ denotes the set of solutions $x:[0,1]\rightarrow M$, such that

\begin{equation} \label{e:generator}
\left\{ \begin{aligned}
         \dot x(t)=X_{H_{L_0, L_1}}(x(t)),\\
         x(0)\in L_0,\hskip 3mm x(1)\in L_1,
                          \end{aligned} \right.
                          \end{equation}

where $H_{L_0, L_1}$ is part of the universal and consistent choice of Floer data for the pair $(L_0, L_1)$.



$\bullet$ \textit{$A_\infty$-compositions:}

We use perturbed holomorphic polygons to define the $A_\infty$-compositions in $\fuk(M)$. Concretely, let

$$\mu^d:hom(L_{d-1},L_d)\otimes\cdots\otimes
hom(L_0,L_1)\rightarrow CF^*(L_0,L_d)$$

for $d\ge 1$ be defined as:

\beq\label{e:mud-noneq}\mu^d(x_d,\cdots,x_1)=
      \sum_{\begin{subarray}{l}
         y\in\cC(L_0,L_d;H_{L_0,L_d})\\ u\in\cM_{\cD_{d,1}}(y;x_d,\cdots,x_1)\\
         \hskip 4mm\text{dim}\cM_{\cD_{d,1}}=0
        \end{subarray}}sign(u)T^{\w(u)}\la y\ra\eeq

which are algebraic counts on rigid holomorphic polygons with $x_1,\cdots,x_d$ as inputs and $y$ as outputs, weighted by the Novikov coefficients determined by the area of the polygon.

Here we recall that $\cM_{\cD_{d,1}}(y;x_d,\cdots,x_1)$ is the solution of \eqref{e:FloerEq} defined by
the universal and consistent choice of Floer data for the $(d+1)$-tuple $(L_0,\cdots, L_d)$,
$sign(u)$ is canonically determined from the orientation of $u$, see \cite{FOOO_Book}\cite[Section 11]{Seidelbook} for a comprehensive account on the sign issues, which we will omit for the rest of the paper. Such compositions $\{\mu^d\}$ satisfies the $A_\infty$-relations as proved again in \cite{FOOO_Book} and \cite{Seidelbook}, which we will not reproduce here.

\subsection{Incorporating the $G$-action}\label{s:Gact}

We now bring in the symplectic $G$-action on $(M,\w)$
and explain term by term our adaption from the previous section
to define the equivariant Fukaya category $\fuk(M)^G$.

$\bullet$ \textit{G-equivariant gradings and spin structures:}

We further restrict ourselves to the class of symplectic manifold $M$ and Lagrangian submanifolds $L$ satisfying the following:


\begin{assumption}\label{a:assumption}\hfill
\begin{enumerate}
\item \label{e:gradingAssumption}$M$ is equipped with a $G$-invariant $N$-Maslov covering for $2\le N\le\infty$.
$L$ is equipped with a $G_L$-equivariant $\Z/N$-grading.


\item\label{e:spinAssumption} We assume the induced
action of $G_L$ on the orthonormal framebundle $O(L)$ lifts to
an action of an associated spin bundle $Spin(L)$.
\end{enumerate}
\end{assumption}

Some discussions on these assumptions seem instructive. The keypoint is that, with these assumptions, any Lagrangian we
consider is a connected component of a lift of Lagrangian brane from the quotient, so that we could easily compare $\fuk(\ov M)$ and $\fuk(M)^G$ using the transfer functor $\cT$ (see Section \ref{s:transfer}).

The assumption on the spin structure is different from that of \cite{Seidelbook}: when we consider $G$-actions on the spin structure, we do not need additional twist between the spin structure restricted on different components, hence is much simpler to deal with. In the language of \cite{CH13}, we use the \textit{spin profile} equals zero. The reader interested in this subject is referred to these nicely written references, which we will not discuss any more.

For the condition (\ref{e:gradingAssumption}), there is a sufficient condition following arguments in \cite{Seidelbook}\cite{CH13}.

\begin{lem}\label{l:grAssm}
 Assume $gcd( ord(G), N )=1$, then $M$ admits an $N$-fold Maslov cover if and only if it admits a $G$-equivariant $N$-fold Maslov cover. An embedded Lagrangian $L$ admits a $\Z/N$-grading if and only if $L$ admits a $G_L$-equivariant $\Z/N$-grading.
\end{lem}

In particular if $N=\infty$, it was shown by Seidel when $G=\Z/2$ and Cho-Hansol for the general case that there is no extra obstruction of obtaining an equivariant grading.

The proof goes as follows. For the Maslov cover part, one could pass to the quotient and conclude that $ord(G)\cdot 2c_1(M/G)=0\in H^2(M/G,\Z/N)$ which is in turn equivalent to $2c_1(M/G)=0\in H^2(M/G,\Z/N)$. Hence one may obtain an $N$-fold Maslov cover from the quotient and pull back to $M$. Alternatively, one may average
the global Maslov class $C_N\in H^1(\sL, \Z/N)$ using the $G$-action
using the invertibility of $ord(G)\in \Z/N$.

For the equivariant grading of $L\subset M$, given an ordinary grading $gr^N:L\rightarrow \sL^N$ supported on the $G$-equivariant Maslov cover, solely its existence implies that $(gr^N)^*(\sL^N)\rightarrow L$ is a trivial covering, or equivalently, is a trivial $\Z/N$-bundle. Assume $g\in G_L$, $gr^N\circ g-gr^N\in\Z/N$ is well-defined and locally constant. Hence if $ord(G)$ is invertible in $\Z/N$, then
$$ord(g)(gr^N\circ g-gr^N)=\sum_{k=1}^{ord(g)}(gr^N\circ g^k-gr^N\circ g^{k-1})=0$$

implies $gr^N\circ g-gr^N=0$, hence the claim.

However, \eqref{e:gradingAssumption} is clearly broader than Lemma \ref{l:grAssm}. Namely, suppose $G=\Z/2$ and $L$ is only $\Z/2$-graded due to its orientability. If $L/G\subset M/G$ is still orientable, it clearly obtains a $\Z/2$-grading which lifts to an equivariant one on $L$.





$\bullet$ \textit{Objects: $G$-Lagrangian branes.}

For any connected embedded Lagrangian brane $L^\#=(L, gr^N, spin_L)$ satisfying Assumption \ref{a:assumption}, we denote the $G$-Lagrangian brane associated
to $L^\#$ as

$$GL^\#:= (\bigcup_{g\in G}g\cdot L, spin_{GL}, gr^N_{GL}), $$

which is the following. Its underlying Lagrangian submanifold is the orbit of $L$ under the $G$-action.  We emphasize that we only consider the underlying set and not the multiplicities caused by non-trivial $G_L$'s. Since the decoration $(gr^N, spin_L)$ is $G_L$-invariant, one then legitimately transfers this invariant grading and spin structure to other components $gL$ using the $G$-action, which consists the decorations $(spin_{GL}, gr^N_{GL})$ on the whole orbit $GL$. Note that since $gL$ and $g'L$
are allowed to intersect when $g\neq g'$, the underlying Lagrangian submanifold is usually immersed.  We will call the orbit elements $gL$ for any $g\in G$ \textit{irreducible components} of $GL^\#$.



$\bullet$ \textit{The morphism groups.}

We first make a heuristic definition of the Floer cochain group
between two $G$-Lagrangian branes as

\beq\label{e:ChainGroup} CF^*(GL_0,GL_1)=\bigoplus_{g_0,g_1\in
G}CF^*(g_0L_0,g_1L_1) \eeq

when there is no isotropy groups for $L_0$ or $L_1$. In the general
case, we take exactly one representative $g_i$ in each coset
elements in $G/G_{L_i}$, $i=0,1$, respectively, in equation
\eqref{e:ChainGroup}.

To make sense of this definition, we need to modify slightly the universal and consistent choice of regular Floer data for the non-equivariant case. First of all, note that we did not include a preferred element as part of the data of the equivariant Lagrangian, so that $GL^\#$ and $G(gL^\#)$ are considered the same object.  Instead, we include this extra piece of information into the Floer data: for each $G$-Lagrangian brane $GL^\#$ we fix one of its irreducible component $L$ and call it the \textit{principal component}.  This choice should be constant over any moduli problem of perturbed holomorphic curves involved in \ref{s:moduli}, and is implicit in the notation $GL$.

The rest of the argument goes as in \cite{Seidelbook}. Given any pair of $G$-equivariant Lagrangian branes $(GL_0^\#, GL_1^\#)$, we take pairs of components $(L_0, gL_1)$ for one element $g$ in each coset $G/G_{L_1}$, where $L_0, L_1$ are the principal components of $(GL_0^\#, GL_1^\#)$. One thus obtains a generic set $\cP^G_{(L_0, gL_1)}\subset \cP^G$ for each $g$ which is regular in the usual sense from Section \ref{s:ET}.  Once a Floer datum is picked for the pair $(L_0,gL_1)$, the $G$-invariance automatically determines the same piece of Floer datum for pairs of the form $(g'L_0,g'gL_1)$, $g'\in G$. Take the intersection of these subsets and denote it as $\cP^G_{GL_0, GL_1}$, any Floer datum in this subset would define \eqref{e:ChainGroup}.
One then inductively extend the Floer data to any $(d+1)$-tuple of $G$-Lagrangian branes as in \cite{Seidelbook} for the non-equivariant case.


$CF^*(GL_0,GL_1)$ admits an obvious
$G$-action preserving degrees from our choices. Now we define the
morphism group as the $G$-invariant part of the Floer cochain group

\beq\label{e:GMor}hom_{\fuk(M)^G}^*(GL_0,GL_1):=CF^*(GL_0,GL_1)^G.\eeq

We emphasize that from the freeness condition, any  morphism in $hom^*(GL_0,GL_1)$ has the form $G\cdot x$ for some
$x\in \bigoplus_{g\in G}CF(L_0,gL_1)$.


$\bullet$ \textit{$A_\infty$-Compositions.}

We continue to use counts of
holomorphic polygons with Hamiltonian perturbation at strip-like
ends to define compositions of $hom^*(GL_0,GL_1)$. The
$G$-invariance of Floer data ensures that the usual composition map
lands on the $G$-invariant part of the Floer cochain group. Namely,
assume a perturbed holomorphic polygon $u\in \cM(x_0;x_1,\cdots,x_d)$
contributes to the composition:

\beq\label{e:mud}\mu^d_G:hom(GL_{d-1},GL_d)\otimes\cdots\otimes
hom(GL_0,GL_1)\rightarrow CF^*(GL_0,GL_d)\eeq

for $x_0\in hom^*(g_0L_0,g_dL_d)$ and $x_i\in
hom^*(g_{i-1}L_{i-1},g_iL_i)$, $i=1,\cdots,d$ for some group elements
$g_i$. Then by the $G$-invariance of our perturbation data, $gu\in
g\cdot\cM(gx_0;gx_1,\cdots, gx_d)$ for any $g\in G$ also contributes
to \eqref{e:mud}.

 To summarize, in the $G$-equivariant case, the composition maps $\{\mu^d_G\}_{d\ge1}$ is no more than an additive enlargment of the usual $A_\infty$-compositions with $G$-invariant regular Floer data. The $A_\infty$ relation for the $\{\mu^d_G\}_{d\ge1}$ hence follows from non-equivariant ones.

$\bullet$ \textit{An algebraic point of view.}

To complete our discussion on the definition of $\fuk(M)^G$ in this section, we would like to relate it to $\fuk(M)$.  There is an obvious functor

$$\wt\iota: \fuk(M)^G\rightarrow Tw(\fuk(M))$$

sending $GL$ to $\oplus_{g\in G} gL$ with vanishing differentials.  Since our choices of Floer data are both $G$-invariant and regular, the Fukaya categories on both sides are defined using the same choice of universal and consistent Floer data, the obvious map on morphism level is also well-defined.  Therefore, by setting $\wt\iota^d\equiv0$ for all $d\ge2$ yields an $A_\infty$-functor.

There is a naive $G$-action on $Tw(\fuk(M))$.  Clearly, $\wt\iota$ reduces to a functor to the fixed category of $Tw\fuk(M)$ as defined in Section \ref{s:algreduction}.  Moreover, this induces an $A_\infty$-functor

$$\iota: Tw(\fuk(M)^G)\rightarrow (Tw\fuk(M))^\fix.$$

$\iota$ is an isomorphism of $A_\infty$-categories: it follows from that each $G$-invariant twisted complex on the right hand side is formed by the $G$-orbit of a direct sum of Lagrangian objects, then the $G$-orbit of each direct summand is an object of $\fuk(M)^G$.  Entries of the differential can also be written similarly as $G$-orbits of elements of Floer cochains.

What we explained above shows the objects of the two categories are one-one correspondent.  The fully faithfulness of $\iota$ can be argued similarly as the differential part.  Therefore, Lemma \ref{l:dfixfixd} shows:

\blem\label{l:eqeq} $D^\pi(\fuk(M)^G)=(D^\pi\fuk(M))^\fix$ is an isomorphism of triangulated categories. \elem


\subsection{The $G$-equivariant generation criterion}

In this section, we give the statement and proof of the $G$-equivariant generation criterion \ref{t:gen} and \ref{t:gen2}.  The split-generation is in the sense of Definition \ref{d:Aidem}, which implies the split generation on derived level.

Our proof follows closely that of \cite{AbGen}.  We start by summarizing several other algebraic operations including coproducts, the open-closed and closed-open string maps in $\fuk(M)^G$ relevant to our proof. The way to define them using various moduli spaces of the form $\cM_{\sR}$ is similar to the $A_\infty$-structure: they send the tensor product of inputs to the tensor products to the outputs, adding the Novikov coefficients corresponding to rigid objects in the moduli space. The algebraic structure of these operations are in turn derived from the degeneration scenarios of relevant moduli spaces.

We will focus on these algebraic operations in the equivariant case.  However, they are exactly the restriction from non-equivariant case to the $G$-invariant part of corresponding (co)chains defined using $G$-invariant Floer data.  Hence we will not reproduce the proof that these maps are (co)chain maps, which can be found in \cite{AbGen}\cite{RS13}\cite{FOOO_Book} \cite{BC1}--in fact \cite{AbGen}\cite{RS13}\cite{Gan13} these maps are even shown to be chain maps for the wrapped cases.  But we will still attach a super(sub)-script to emphasize when the $G$-equivariant cases is being considered.

\nin\textbf{The coproduct $\Delta_G$:} Coproducts are defined by the moduli problem of $\cM_{\cD_{n,2}}$. Explicitly, this is a degree-$n$ homomorphism of $A_\infty$-bimodules

$$\Delta_G:G\cB\rightarrow \cY_{GK}^l\otimes \cY_{GK}^r$$

\nin for any full subcategory of $G\cB\subset \fuk(M)^G$ and $GK\in Ob\fuk(M)^G$. By definition, this consists of a collection of maps:

\begin{multline}\label{e:Coprod}\Delta_G^{r|1|s}:CF^*(GL_{r-1},GL_r)
\otimes\cdots\otimes CF^*(GL_0,GL_1)\otimes CF^*(GL_{|0},GL_0)\\
\otimes CF^*(GL_{|1},GL_{|0})\otimes\cdots\otimes CF^*(GL|_{|s},GL_{|s-1})\rightarrow CF^*(GK,GL_r)\otimes CF^*(GL_{|s},GK)\end{multline}

\nin for all involved Lagrangians $GL_\gamma\in G\cB$, and satisfies the certain $A_\infty$ equation (see \cite[Section 4.2]{AbGen})


\nin$\bullet$\hskip 2mm\textbf{The open-closed string map $\OC^G$:} Let $CF^*(M)^G$ denote the invariant part of $CF^*(M)$.  The open-closed
string map is defined by the moduli problem of $\cM_{\cD_{n,0}^-}$. It gives chain level homomorphisms:

\begin{equation} \label{eq:open_closed_length_d} \OC^G_{d} \co CW^{*}(GL_{d-1},GL_0) \otimes \cdots \otimes CF^{*}(GL_1, GL_2) \otimes CW^{*}(GL_0, GL_1) \to CF^{*}(M)^G \end{equation}
which shift degree by $n-d+1$ and are the components of a degree $n$ chain map
\begin{equation}
\OC^G \co CC^G_{*}(G\cB, G\cB) \to CF^{*}(M)^G.
\end{equation}

Here the left hand-side is the cyclic bar complex of $G\cB$ equipped with the differential computing Hochschild homology of $\fuk(M)^G$.

\nin$\bullet$\hskip 2mm\textbf{The closed-open string map $\CO$:}

The moduli problem of $\cM_{\cD_{0,1}^+}$ defines the closed-open string map, which is a chain map between Floer cochain complexes:

\begin{equation}
  \label{eq:closed_to_open}
\CO \co CF^*(M)^G \to CF^{*}(GK,GK)
\end{equation}

\nin where $GK$ is the unique Lagrangian label on $\partial S$ for $S\in\cD_1^+$.

\nin$\bullet$\hskip 2mm\textbf{The homotopy $\cH$:}

The last operation $\cH$ is a chain map

$$\cH:CC_*(G\cB, G\cB)\rightarrow CF^*(K)^G[n]$$

\nin between the cyclic bar complex and the Floer complex defined by the moduli problem of $\cM_{\ov\cC_n^-}$, while inputs on the left hand side comes from Lagrangian labels on $\{|z|=r>1\}\subset \partial S$, $S\subset \cC_n^-$ and $K$ is the unique Lagrangian label on $\{|z|=1\}$.

Next we define a purely algebraic morphism. Recall that given $\cL, \cR$ a left (resp. right) $G\cB$-module, then the tensor product over $G\cB$ is a chain complex:

\begin{equation}\label{e:tensor}\cR\otimes_{G\cB}\cL=
\bigoplus_{GL_0,\cdots, GL_d\in Ob(G\cB)}\cR(L_d)\otimes CF^*(L_{d-1},L_d)\otimes\cdots\otimes CF^*(L_0,L_1)\otimes\cL(L_0)
\end{equation}

\nin with differential

\begin{multline*}
  p \otimes a_d \otimes \ldots \otimes a_1 \otimes q \mapsto \sum p \otimes a_{d} \otimes \cdots \otimes a_{\ell+1} \otimes \mu^{\ell|1}(a_{\ell}, \ldots, a_{1}, q) \\
+ \sum (-1)^{\deg(q) + \maltese_{1}^{\ell}} \mu^{1|d-\ell} (p , a_{d}, \ldots, a_{\ell+1}) \otimes a_{\ell} \otimes \cdots \otimes a_{1} \otimes q \\
+ \sum (-1)^{\deg(q) + \maltese_{1}^{\ell}}p \otimes a_{d} \otimes \cdots \otimes a_{\ell+k+1} \otimes \mu^{k}(a_{\ell+k}, \ldots, a_{\ell+1}) \otimes a_{\ell} \otimes \cdots \otimes a_{1} \otimes q.
\end{multline*}

The bimodule morphism $\Delta_G$ induces at the level of Hochschild chains a homomorphism $CC_{*}(\Delta_G): CC_*^G(G\cB)\rightarrow \cY^r_{GK}\otimes_{G\cB}\cY^l_{GK}$:

\begin{multline}
  \label{eq:hochschild_map_induced_by_coproduct}
CC_{*}(\Delta_G) (a_{d} \otimes \ldots \otimes a_{1}) = \\ \sum (-1)^{\diamond} \cI\left( \Delta_G^{r|1|s}(a_{r}, \ldots, a_{1}, \underline{a}_{d}, a_{d-1}, \ldots, a_{d-s}) \otimes a_{d-s-1} \otimes \cdots \otimes a_{r+1} \right)
\end{multline}
where $\cI$ is the maps which reorders the factors
\begin{equation*}
\cI( q \otimes p \otimes a_{d-s-1} \otimes \cdots \otimes a_{r+1}) = (-1)^{\circ} p \otimes a_{d-s-1} \otimes \cdots \otimes a_{r+1} \otimes q
\end{equation*}
and the signs are given by the formulae
\begin{align}
\label{eq:sign_delta_induce_cyclic} \diamond & = \maltese_{1}^{r} \cdot (1+ \maltese_{r+1}^{d}) + n \maltese_{r+1}^{d-s-1} \\
  \label{eq:sign_reorder} \circ & = \deg(q) ( \deg(p) + \maltese_{r+1}^{d-s-1})\\
  \maltese_s^t &=\sum_{s\le j\le t}||a_j||.
\end{align}

We define $HH_{*}(\Delta_G) $ to be the map induced by $ CC_{*}(\Delta_G) $ on homology groups.

The last bit of information we need is the composition map,
which is a chain map of degree $0$:

\begin{align}
  \label{eq:composition_left_right}
\mu \co \cY^{r}_{GK} \otimes_{G\cB} \cY^{l}_{GK} & \to CF_G^{*}(GK,GK) \\
p \otimes a_{d} \otimes \cdots \otimes a_{1} \otimes q & \mapsto (-1)^{\deg(q) + \maltese_{1}^{d}} \mu_G^{d+2}(p, a_{d}, \ldots, a_{1}, q).
\end{align}

What is important to us is the following proposition.

\bprop\label{p:CommDiag}
The following diagram commutes up to homotopy and a sign $(-1)^{n(n+1/2)}$:
\beq\label{e:commuDiag}
\xymatrix{
CC_*^G(G\cB) \ar[d]^{\OC^G} \ar[rr]^{CC_*(\Delta_G)} && \cY^r_{GK}\otimes_{G\cdot\cB}\cY^l_{GK} \ar[d]^{\mu} \\
CF^*(M)^G \ar[rr]^{\CO^G} && CF^*(GK,GK)^G } \eeq

\nin Moreover, the chain homotopy is defined by $\cH$.
\eprop

\bpf For the non-equivariant case, this is proved in \cite[Proposition 1.3]{AbGen} and \cite[Section 8.6]{RS13}. But in our case it is actually much simpler: consider the Gromov bordification of dimension $1$ moduli of type $\cM_{\cC^-_d}$. Its boundaries consist of the following four types of moduli spaces:

\begin{align}
\label{e:bdry1}&\cM_{\cD_{n,0}^-}\times\cM_{\cD^+_{0,1}} \\
\label{e:bdry2}&\cM_{\cC_{d_1}^-}\times\cM_{\cD_{d_2,1}}, \hskip 3mm d_1+d_2=d\\
\label{e:bdry3}&\cM_{\cC_{d}^-}\times\cM_{\cD_{1,1}}\\
\label{e:bdry4}&\cM_{\cD_{d_1,2}^-}\times\cM_{\cD_{d_1,1}}, \hskip 3mm d_1+d_2-2=d.
\end{align}

Here \eqref{e:bdry2},\eqref{e:bdry3} are responsible for the chain homotopy part

$$(-1)^n\mu_{\fuk(M)^G}^1\circ\cH+\cH\circ b,$$

\nin where $b$ is the differential of the Hochschild chain complex and \eqref{e:bdry1},\eqref{e:bdry4} accounts for the difference

$$\mu\circ CC_*(\Delta_G)-\CO^G\circ\OC^G.$$

Hence we obtain

\beq\label{e:homotopyEq} (-1)^n\mu^1\circ\cH+\cH\circ b+\mu\circ CC_*(\Delta_G)-\cC\cO^G\circ\cO\cC^G=0\eeq

This is exactly equation (6.9) in \cite{AbGen}. Our simplification relies on the fact that there is
no rescaling trick necessary for our applications, hence the moduli space counting $\mu$ and $CC_*(\Delta)$
can be glued directly. This process is essentially the homotopy $\cH_2$ in \cite{AbGen}.

\epf

From Proposition \ref{p:CommDiag}, we have a commutative diagram up to a sign $(-1)^{\frac{n(n+1)}{2}}$:

\beq\label{e:commuDiag}
\xymatrix{
\HHGGB \ar[d]^{H^*(\OC^G)} \ar[rr]^{HH_*(\Delta_G)} && H^*(\cY^r_{GK}\otimes_{G\cdot\cB}\cY^l_{GK}) \ar[d]^{H^*(\mu)} \\
HF^*(M)^G \ar[rr]^{H^*(\CO^G)} && HF^*(GK,GK)^G } \eeq

We are now ready to explain the proof of Theorem \ref{t:gen}.  Note first $id\in HF^*(M)$ is automatically $G$-invariant as the fundamental class. Assume that $H^*(\OC^G)$ hits the identity $id\in HF^*(M)^G$.  Given an embedded Lagrangian $K\subset M$, $H^*(\CO)$ sends $id_{HF^*(M)}$ to $id_{HF^*(gK)}$ (see for example \cite{SeidelBias}) for each irreducible component $gK$.  Therefore, in the equivariant setting, $H^*(\CO^G)(id_{HF^*(M)^G})=id_{HF^*(GK)}=\sum_{g\in G}id_{HF^*(gK)}$ from applying the $G$-action.



The following lemma is due to Abouzaid \cite[Lemma 1.4]{AbGen}, see also \cite[Proposition 2.6]{Gan13}.

\blem\label{l:GenLem} Given an $A_\infty$ full subcategory $\cA'$ in $\cA$, and $K\in Ob(\cA)$. Then $H^*(\mu):\cY_K^r\otimes_{\cA'}\cY_K^l\rightarrow Hom(K,K)$ hits the identity iff $\cA'$ split generates $K$.

\elem

In our $G$-invariant context, the lemma still applies without any modifications: we only regard $\cF uk(M)^G$ as an ordinary $A_\infty$-category
and $G\cB$ generates a full subcategory in the ordinary sense (using morphisms in $\cF uk(M)^G$, which are by definition $G$-invariant). Hence, we have proved:

\bthm\label{t:gen2} Let $G\cB$ be a $A_\infty$ full subcategory of the $\cF uk(M)^G$. Assume that $$\xymatrix{HH_*^G(G\cB) \ar[rr]^{H^*(\OC_G)}&&HF^*(M)^G}$$
has its image containing $id_{HF^*(M)^G}$, then $G\cB$ split-generates $\cF uk(M)^G$.
\ethm

\bcor\label{c:gen} Suppose $\cB\subset Ob(\fuk(M))$ split generates $\fuk(M)$, then $G\cB$ split generates $\fuk(M)^G$.
\ecor

\bpf From the assumption and Lemma \ref{l:GenLem}, we have $H^*(\OC):HH_*(\cB)\rightarrow HF^*(M)$ hits the identity. Take a chain representative of a preimage
$\alpha\in CC_*(\cB)$ so that $\OC(\alpha)=e$ and $[e]=id_{HF^*(M)}$. Then clearly $G\alpha\in CC_*(G\cB)$ and $[\OC^G(G\alpha)]=[Ge]=ord(G)\cdot id_{HF^*(M)}$, which shows $G\cB$ split generates $\fuk(M)^G$ when $|G|$ is invertible.

\epf


\subsection{The transfer functor}\label{s:transfer}

The goal of this section is to define a transfer functor $\cT$ and relate the equivariant Fukaya category to the ordinary Fukaya category of the quotient.  Our main result reads:

\begin{thm}\label{t:functor2} There is an $A_\infty$ functor
\beq\label{e:Tfunctor}\cT:\fuk(\ov M)\rightarrow\fuk(M)^G.\eeq

\noindent which is full and faithful. In particular, if there is a subcollection $\ov{\cB}\subset Ob(\fuk(\ov M))$ so that $\cT(\ov \cB)$
resolves the diagonal, then $\ov\cB$ split generates $\fuk(\ov M)$.

We will call $\cT$ the \textbf{transfer functor}.

\ethm

We now explain the definition of the transfer functor $\cT$.

\textbf{Object level:} For $\ov L\in Ob(\fuk(\ov M))$, we assume it comes with a chosen spin structure $spin(\ov L)$ and grading $gr^N_{\ol L}$.
Let $L$ be a component of $\pi^{-1}(\ov L)$. $spin(\ov L)$ naturally lifts to a $G$-invariant spin structure on $\bigcup_{g\in G}gL$.

The same story holds for the grading: the Grassmannian bundle $Gr(T\ov M)|_{\ov L}$ is lifted to $Gr(TM)|_L$ by $\pi$, hence the section
$gr^N_{\ol L}$ is lifted equivariantly to one on $Gr(TM)|_L$ as in Lemma \ref{l:grAssm}. A possible confusion is that such lifting seems to give a possibility of improving the grading (for example, the vanishing of mod-$N$ global Maslov class might be improved to its vanishing in $\Z$-coefficient after the lift), but this cannot be achieved in a $G$-equivariant way unless the improved
grading can already be realized before the lift.

\textbf{Morphisms:} We require $\cT^d\equiv0$ for all $d\geq2$. $\cT^1$ is defined by lifting $\bar z\in hom(\ov L_0, \ov L_1)$ to $Gz\in hom(GL_0,GL_1)$, i.e. $\cT^1(\bar z)=Gz$. To verify that $\cT$ indeed defines an $A_\infty$ functor, we want to see that

\beq \label{e:TransferEq}\cT^1(\mu_{\ol M}^d(\bar z_d,\cdots, \bar z_1))=\mu^d_G(Gz_d,\cdots Gz_1)
\eeq

\noindent holds for any $\bar z_i\in CF^*(\ov L_{i-1}, \ol L_i), i=1, \cdots, d$. Since both sides of the equality are $G$-invariant, it suffices to compare the coefficient of a fixed lift $z_0$ of some fixed output $\bar z_0$ of $\mu^d_{\ov M}$. Notice the following correspondence:
on the one hand, given any holomorphic polygon in $\bar u\in\cM(\bar z_0; \bar z_d,\cdots,\bar z_1)$, it has precisely $|G|$ lifts which are still holomorphic polygons in $M$ with respect to the lifted data. Among them there is a unique polygon $u$ contributing to $\cM(z_0; z_d,\cdots,z_1)$ for some lifts $z_1,\cdots,z_d\in M$ of $\bar z_1,\cdots, \bar z_d\in\ol M$.  The virtual dimension, area and sign of $u$ are all the same as $\bar u$ due to the freeness assumption. On the other hand, given a holomorphic polygon $u\in \cM(z_0; g_dz_d,\cdots,g_1z_1)$ with $G$-equivariant Floer data for some $g_1,\cdots,g_d\in G$, it clearly descends to $\bar u\in\cM(\bar z_0; \bar z_d,\cdots,\bar z_1)$, also with the same virtual dimension and Novikov coefficients on up- and down-stairs as argued above.
To summarize, we saw that the correspondence from $u$ to $\bar u$ is $ord(G)$-to-one, and singling out the specific lift $z_0$ upstairs gives equality of Novikov coefficients
$$\langle\bar z_0; \mu^d(\bar z_d,\cdots,\bar z_1)\rangle=\sum_{g_i\in G}\langle z_0; \mu^d(g_dz_d,\cdots,g_1z_1)\rangle=\ip{z_0;\mu_M^d(Gz_d,\cdots,Gz_1)}.$$

Summing the above equality over its $G$-orbits for all outputs of the form $gz_0$, $g\in G$ yields the desired equality \eqref{e:TransferEq}.

To this end, we have shown that $\cT$ defines an $A_\infty$ functor, the fact that $\cT$ is fully faithful is trivial since (1) embedded Lagrangians always lifts to embedded Lagrangian submanifolds due to the freeness of $G$-action and, (2) $\cT^1$ sends a basis to basis on the morphism level since all morphisms in $\fuk(M)^G$ has the form $G\cdot x$. This concludes Theorem \ref{t:functor2}. \qed

\bcor\label{c:downstairsGen} If $\cB\in Ob(\fuk(M))$ resolves the diagonal, and $\pi(L)$ is embedded for all $L\in\cB$, then the collection
$\pi(\cB)$ split generates $\fuk(\ov M)$.  Hence $D^\pi(\fuk(\ol M))\cong D^\pi(\fuk(M)^G)\cong (D^\pi(\fuk(M)))^\fix$ from Lemma \ref{l:eqeq}.
\ecor

\bpf We already have a full and faithful functor $\cT:\fuk(\ov M)\rightarrow \fuk(M)^G$. The collection $\cT(\pi(\cB))=G\cB$ by definition. From Corollary \ref{c:gen} this collection generates the whole $\fuk(M)^G$, hence the claim.

\epf

\section{Applications to homological mirror symmetry}

\subsection{Special isogenous tori}\label{s:sympSI}

\subsubsection{Special isogenous symplectic tori}

Let $\bar\alpha=(\alpha_1,\cdots,\alpha_n)\in \R^n$, we will denote
the split torus $T(\bar\alpha):=T(\alpha_1)\times\cdots\times
T(\alpha_n)$, where $T(\alpha_i)$ is the symplectic $2$-torus with
area $\alpha_i$. For the purpose of being explicit, we parametrize
any symplectic $2$-torus with area $A$ as $\{(s,t)\in S^1(A)\times
S^1(1)\}$. We also tacitly rearrange coordinates of $T(\bar\alpha)$
so that the $s$-coordinates are aligned at the first $n$ components,
i.e. $T(\bar\alpha)=(S^1(\alpha_1)\times
S^1(\alpha_2)\times\cdots\times S^1(\alpha_n)\times S^1(1)^n,\sum
ds_i\wedge dt_i)$.

\bdf\label{d:symptori} A \textit{special automorphism} is a finite
$\Z/l$-action for some $l\in\Z$ on a split symplectic torus $T(\bar\alpha)$ with
generator

$$g:T(\alpha)\rightarrow T(\alpha),\hskip 3mm g(s_1,t_1,\cdots,s_n,t_n)\rightarrow
(s_1+\frac{\alpha_1}{l},t_1,\cdots,s_n+\frac{\alpha_n}{l},t_n).$$

We will denote $T(\bar\alpha)_l:=T(\bar\alpha)/(\Z/l)$ and call
it a \textit{special isogenous torus with one factor}. A \textit{special isogenous torus} is then
a finite product of such tori.

\edf

Note that a split torus is also a special isogenous one with $l=1$. Let us now focus on the special isogenous tori with one factor at the moment.  $T(\bar\alpha)_{l}$ can be realized as a lattice quotient of
$(\R^{2n},\w_{std})/\Gamma^{\bar\alpha}_{l}$. Explicitly,

$$\Gamma^{\bar\alpha}_{l}=\Z e_1+\cdots+\Z e_{2n},$$

where

\beq\label{e:basis}\begin{aligned}
&e_{1}=(\alpha_1/l,\cdots,\alpha_n/l,0,\cdots,0)^T\\
&e_i=(0,\cdots,0,s_i=\alpha_i,0,\cdots,0)^T,\hskip 3mm 2\le i\le n,\\
                  &e_{n+j}=(0,\cdots,0,t_j=1,0,\cdots,0)^T,\hskip 3mm 1\le j\le n,                  \end{aligned}\eeq

This lattice contains the split lattice formed by the one formed by replacing $e_1$ by $(\alpha_1,0,\cdots,0)^T$ as an index-$l$ subgroup.

Now take a linear transformation defined by the block matrix

$$\wt M_{l}^{\bar\alpha}=\left( \begin{array}{c:c} M_{l}^{\bar\alpha} & 0 \\ \hdashline
0 & I_n \end{array} \right), $$

where

\[
M_{l}^{\bar\alpha}=\left(
 \begin{array}{cccc}
   \frac{l}{\alpha_1} \\
   -\frac{1}{\alpha_1} & \frac{1}{\alpha_2} & \text{\huge0}& \\
   \vdots & &\cdots\\
   -\frac{1}{\alpha_1}&\text{\Large 0} & & \frac{1}{\alpha_n}\\

 \end{array}
\right).
\]

which sends $\Gla$ to the standard lattice $\Z^{2n}\subset
\R^{2n}$. The pull-back symplectic form on $\R^{2n}/\Z^{2n}$,
regarded as a bilinear form on $\R^{2n}$ due to translation
invariance, can then be represented by a matrix

$$\Omega_{\bar l}^{\bar\alpha}=(\wt M_{l}^{\bar\alpha})^T\left( \begin{array}{cc} 0 & I_n \\
-I_n & 0 \end{array} \right)\wt M_{l}^{\bar\alpha}=\left(\begin{array}{cc} 0 & (M_{l}^{\bar\alpha})^T \\
-M_{l}^{\bar\alpha} & 0 \end{array} \right) .$$

Hence, if $T(\bar\alpha)_{l}$ and $T(\bar\alpha')_{l'}$
are symplectomorphic, the two corresponding matrices
$\Omega_{l}^{\bar\alpha}$ and $\Omega_{l'}^{\bar\alpha'}$
are in the same $GL(2n,\Z)$-congruence class (recall that we already argued in the introduction that linear symplectic tori are symplectomorphic only if they are symplectomorphic by linear transformations).  The converse of the above argument works equally well.  To summarize:

\blem\label{l:LSeqCong} $T(\bar\alpha)_{l}$ and $T(\bar\alpha')_{l'}$
are symplectomorphic if and only if  the two corresponding matrices
$\Omega_{l}^{\bar\alpha}$ and $\Omega_{l'}^{\bar\alpha'}$
are in the same $GL(2n,\Z)$-congruence class.

\elem

Finally, note that the results above apply to a general special isogenous torus--one simply replaces the matrix $M_l^{\bar \alpha}$ by a block matrix with blocks of the same shape.

\bex
 There is an obvious generalization of the special automorphism construction.  Take
 $\bar l=(l_1,\cdots,l_n)\in\Z^n$ with $gcd(l_1,\cdots,l_n)=1$ and $l=lcm(l_1,\cdots,l_n)$.  Define a $\Z/l$-action by the generator

 $$g:T(\alpha)\rightarrow T(\alpha),\hskip 3mm g(s_1,t_1,\cdots,s_n,t_n)\rightarrow
(s_1+\frac{\alpha_1}{l_1},t_1,\cdots,s_n+\frac{\alpha_n}{l_n},t_n),$$

\hskip3mm We will show that such a construction does not provide any new examples.  Let us denote the result of such a quotient as $T(\bar\alpha)_{\bar l}$.  $T(\bar\alpha)_{\bar l}$ can be realized as a lattice quotient of
$(\R^{2n},\w_{std})/\Gamma^{\bar\alpha}_{\bar l}$. Explicitly,

$$\Gamma^{\bar\alpha}_{\bar l}=\Z e_1+\cdots+\Z e_{2n}+\Z e_{2n+1},$$

where

\beq\label{e:basis}\begin{aligned}
&e_{1}=(\alpha_1/l_1,\cdots,\alpha_n/l_n,0,\cdots,0)^T\\
&e_i=(0,\cdots,0,s_i=\alpha_i,0,\cdots,0)^T,\hskip 3mm 2\le i\le n,\\
                  &e_{n+j}=(0,\cdots,0,t_j=1,0,\cdots,0)^T,\hskip 3mm 1\le j\le n,\\
                  &e_{2n+1}=(\alpha_1,0,\cdots,0)^T.
                  \end{aligned}\eeq

\hskip3mm Of course, this is not an integral basis of the lattice. One may make a simplification of $\Gamma^{\bar\alpha}_{\bar l}$ as follows. First, for given $i\ge1$, multiplying $lcm(l_1,\cdots,\hat l_i,\cdots,l_n)$ to $e_1$ gives a lattice point $(*,\alpha_i\cdot lcm(l_1,\cdots,\hat l_i,\cdots,l_n)/l_i,*)^T$, where $*$ denotes entries of integral multiples of $\alpha_j$, $j\ne i$. This implies that $(0,\cdots,0,\alpha_i/\wt l_i,0,\cdots,0)^T\in\Gamma^{\bar\alpha}_{\bar l}$ for $\wt l_i=l_i/l_i'$, $l_i'=gcd(l_i,lcm(l_1,\cdots,\hat l_i,\cdots,l_n))$. Let $\alpha_i'=\alpha_i/\wt l_i$, we have proved $\Gamma^{\bar\alpha}_{\bar l}$ is spanned by

$$\begin{aligned} &e'_i=(0,\cdots,0,s_i=\alpha_i',0,\cdots,0)^T,\hskip 3mm 1\le i\le n,\\
                  &e_{n+j}=(0,\cdots,0,t_j=1,0,\cdots,0)^T,\hskip 3mm 1\le j\le n, \\
                  &e_{2n+1}=(\alpha_1'/l'_1,\cdots,\alpha'_n/l'_n,0,\cdots,0)^T.
                  \end{aligned}$$

If $l_k'=1$ for certain $1\le k\le n$, one may then reduce $e_{2n+1}$ by setting $s_k=0$ and repeat the same procedure to the rest of the coordinates of $e_{2n+1}$. Note that $l_i'\le l_i$ and the strict inequality holds for at least one $i$ unless $l_i=l_j$ holds for all $i,j$ pairs.  Therefore, after finite steps one arrives at an integral basis of $\Gamma^{\bar\alpha}_{\bar l}$ in the form of

 $$\begin{aligned}
&e''_{1}=(\alpha_1''/l^\circ,\cdots,\alpha''_m/l^\circ,\alpha''_{m+1},\cdots,\alpha''_{n},\cdots,0)^T, \hskip 2mm 1\le m\le n, l^\circ\in\Z.\\
 &e''_i=(0,\cdots,0,s_i=\alpha_i'',0,\cdots,0)^T,\hskip 3mm 2\le i\le n,\\
                  &e''_{n+j}=(0,\cdots,0,t_j=1,0,\cdots,0)^T,\hskip 3mm 1\le j\le n, \end{aligned}$$

 \hskip3mm This means $T(\bar\alpha)_{\bar l}$ is always the product of at most two special isogenous tori $T_1\times T_2$ as in Definition \ref{d:symptori}, while one of them is a split torus.  
 \qed
\eex

\subsubsection{Special isogenous analytic tori}\label{s:analSI}

To explain the special isogenous tori in the analytic category, we
will need to recall some notions from rigid analytic geometry
invented by Tate \cite{Ta71}. Readers wishing for a more thorough treatment on this well-studied subject may consult \cite{FP04} for example. All
varieties involved in this section are assumed to be over the
universal Novikov field with complex coefficients $\Lambda=\Lambda_{\C}$ .

Recall first that $\Lambda$ is endowed with a valuation

$$\sigma(\sum_i a_iT^{\lambda_i})=\lambda_0.$$

We interchangeably refer to a \textit{split algebraic torus of
dimension $n$} or its analytification as
$T^{an}=(\Lambda^*)^n=Spec\Lambda[z_1,z_1^{-1},\cdots,z_n,z_n^{-1}]$.
A \textit{character} on $\Tan$ is a algebraic group homomorphism
$\chi:\Tan\rightarrow \Lambda^*$, which has the form $z^{\bf
k}=z_1^{k_1}\cdot z_2^{k_2}\cdots z_n^{k_n}$ for ${\bf l}=(l_1,\cdots,l_n)\in\Z^n$. We denote the set of characters as $X(\Tan)$.

The valuation extends to $\Tan$ as a map

$$\sigma:\Tan\rightarrow\R^n.$$

A \textit{lattice} $\Gamma\subset\Tan$ is a subgroup such that $\sigma:\Gamma\rightarrow\R^n$ is injective and its image forms a lattice of full rank in the classical sense. By abuse of notation, we sometimes use the notation $\Gamma$ (or lattice points in $\Gamma$) and $\sigma(\Gamma)$ interchangeably when the context is clear. One then considers the quotient $\sT=\Tan/\Gamma$. One has

\bprop[\cite{FP04},6.4.1]\label{p:LQisScheme} $\sT$ is a separated and proper $\Lambda$-scheme.
\eprop

Indeed $\sT$ has the additional structure of a \textit{rigid analytic space}. $\sT$ thus defined will be called an \textit{analytic torus} or simply a \textit{lattice quotient}.

As in the case over $\C$, not all lattice quotients defines an abelian variety. A classical sufficient and necessary condition is given by the Riemann matrix condition (see for example \cite{BL04}). In the rigid analytic setting, we have the following:

\bthm[\cite{FP04},6.6.1]\label{t:LQisAb} $\sT$ is an abelian variety if and only if there is a homomorphism $\vp:\Lambda\rightarrow X(T)$, such that the bilinear form

$$\ip{\lambda,\lambda'}:=\sigma(\vp(\lambda)(\lambda')),\hskip3mm \lambda,\lambda'\in\Gamma$$

is symmetric and positive definite.

\ethm

\begin{exampledefinition}\label{ex:si} Let the lattice $\Gla=\la V_1=(q^{\alpha_1/l},q^{\alpha_2/l},\cdots,q^{\alpha_n/l}),V_2=(1,q^{\alpha_2},1,\cdots,1)\cdots,V_n=(1,\cdots,1,q^{\alpha_n})\rangle$, where $\bar\alpha=(\alpha_1,\cdots,\alpha_n)\in\R^n$. We may define


\begin{align*}
&\vp(V_1)=z_1\cdot z_2\cdot\cdots\cdot z_n,\\
&\vp(V_i)=z_i^l,\hskip 3mm \textit{for }i\ge2.
\end{align*}

Given any lattice points
\begin{align*}
&\lambda_{m}=(q^{\frac{m_1\alpha_1}{l}},q^{\frac{m_1\alpha_2}{l}+m_2\alpha_2},\cdots,q^{\frac{m_1\alpha_n}{l}+m_n\alpha_n})
&\lambda_{m'}=(q^{\frac{m'_1\alpha_1}{l}},q^{\frac{m'_1\alpha_2}{l}+m'_2\alpha_2},\cdots,q^{\frac{m'_1\alpha_n}{l}+m'_n\alpha_n}),
\end{align*}

one may compute

\begin{align*}
\ip{\lambda_m,\lambda_{m'}}&=\sigma((z_1^{m_1}\prod_{i=2}^{n}z_i^{m_il+m_1})(\lambda_{m'}))\\
&=\sum_{i=2}^n(m_il+m_1)(m_1'\alpha_i/l+m_i'\alpha_i)+\frac{\alpha_1}{l_1}\cdot m_1'\cdot m_1\\
&=\sum_{i=2}^n\frac{\alpha_i}{l}(m_il+m_1)(m_i'l_i+m_1')+\frac{\alpha_1}{l}\cdot m_1'\cdot m_1
\end{align*}

which is clearly a positive definite quadratic form. This verifies $(\Lambda^*)^n/\Gla$ is an abelian variety over $\Lambda$.
We will denote this lattice quotient as $A(\bar\alpha)_{l}$ and call a finite product of these tori a \textit{special isogenous analytic torus}.  When $l=1$, we will simply refer to $A(\bar\alpha)=A(\alpha_1)\times\cdots\times A(\alpha_n)=\Lambda^*/\ip{T^{\alpha_1}}\times\cdots\times\Lambda^*/\ip{T^{\alpha_n}}$
as a \textit{split analytic torus}.

Note that a special isogenous analytic torus is indeed isogenous to a product torus $A(\bar\alpha)_l$ by considering the following free action of $\Z/l$ generated by

$$g^{an}:A(\bar\alpha)\rightarrow A(\bar\alpha),\hskip 3mm g(\lambda_1,\cdots,\lambda_n)\rightarrow
(\lambda_1\cdot T^{\alpha_1/l},\cdots,\lambda_n\cdot T^{\alpha_n/l}).$$

The verification is straightforward hence omitted. We will see that this is the mirror group action of the symplectic side (Definition \ref{d:symptori}).

\end{exampledefinition}

\subsection{Homological mirror symmetry of product tori}\label{s:MSPT}


We briefly recall the homological mirror symmetry of a product torus following \cite{AS10}.  We also borrowed expositions in the case of a two-torus from \cite{Ha13} and Auroux's lecture notes \cite{AuNotes}.


Given a torus $T(\alpha)$ with area $\alpha$, parametrized by $(s,t)\in\R/\alpha\Z\times\R/\Z$. We denote longitudes and meridians as:

$$L_{\mathbbm{m},a}=\{s=a\},\hskip 3mm L_{\mathbbm{l},b}=\{t=b\},$$

\noindent respectively. It is shown in \cite{PZ98, AS10} etc. that $(T(\alpha),\w_{std})$ has the Tate curve
$A(\alpha)$ as the mirror, which can be identified with its analytification $\Lambda_\C^*/\langle T^\alpha\rangle$. Note that any $\zeta\in\Lambda_\C^*/\ip{T^\alpha}$
can be written as $\zeta=T^{\sigma(\zeta)}\cdot(a_0+a_1T^{s_1}+\cdots)$, where we take the $mod-\alpha$ reduction $\sigma(\zeta)\in\R/\alpha\Z$ and $s_i>0$ for all $i\ge1$.




We describe the mirror functor at least on the object level for a set of generators. Namely, on the symplectic side, \cite{AS10} showed that $\{L_{\bl,0}, L_{\bm,0}\}$ split generates $D^{\pi}\fuk(T(\alpha))$.
 To describe objects in the derived category as geometric objects, one must incorporate an additional piece of data of \textit{local systems}, which is more or less a routine procedure nowadays, but we will avoid introducing more notations here and only consider them as objects introduced by taking the Karoubi completion. Clearly, the split generation holds true also for any pair of longitude and meridian $\{L_{\bl,b}, L_{\bm,a}\}$ by applying an appropriate symplectomorphism.

On the analytic side, it is known (see \cite{Or09}) that $\cO$ and $\cO_p$ split generates the bounded derived category of coherent sheaves $D^bCoh(A(\alpha))$ for any $p\in A(\alpha)$.
 Moreover, the mirror functor constructed in \cite{PZ98} can be used to describe mirrors of longitudes and meridians. Denote $p_b=e^{2\pi b\sqrt{-1}}\in\C^*\subset\Lambda^*$, one has the following correspondence via the mirror functor:

\beq\label{e:mirror}\begin{array}{l}
 L_{\bl,b}\longleftrightarrow \cO(p_b-p_0),\\
 L_{\bm,a}\longleftrightarrow \cO_{T^a\cdot p_0}.\end{array}\eeq

In short, the longitudes correspond to a degree zero line bundle; while the meridians correspond to skyscraper sheaves of closed points.

Abouzaid and Smith \cite[Section 7,8]{AS10} generalizes this split generation to any product tori.  For $T(\bar\alpha)=\prod_{i=1}^n T(\alpha_i)$, the split generation holds for the set of product Lagrangians of the form

$$L^{W,t}:=L_{w_1,t_1}\times\cdots\times L_{w_n,t_n}.$$

Here $W=w_1\cdots w_n$ is a word of length $n$ consisting of $\bm$ and $\bl$; $t=(t_1,\cdots,t_n)$, where $t_i$ is a given sequence of arbitrary real numbers within its own range of values.   This corresponds to a split generation result of $D^bCoh(A(\alpha))$ (again using \cite[Theorem 4]{Or09}, see \cite[Lemma 8.1]{AS10}), where $L^{W,t}$ corresponds to the coherent sheaf

$$\cE^V:=\cE_1\boxtimes\cdots\boxtimes\cE_n,$$

where $\cE_i$ is the mirror of $L_{w_i,t_i}$ under the correspondence \eqref{e:mirror}.
The upshot is that, given any finite set of Lagrangians of the form $L^{W,t}$ without repetition in words $W$'s, and consider the subcategory $\cA$ consisting of twisted complexes they split generate in $D^\pi\fuk(T(\bar\alpha))$, there is a fully faithful functor

\beq\label{e:Afunctor} \cF: \cA\rightarrow D^b(A(\bar\alpha)).\eeq

Also, the following homological level assertion will be useful.

\blem\label{l:isomorphic}$\cA$ is isomorphic to its image.\elem

   This is obvious from definition: the images of $L^{W,t}$ are pairwisely distinct objects for different value of $W$ and $\cF$ is fully faithful.  Now $\cA$ is equivalent to $D^\pi Fuk(T(\bar\alpha))$ if and only if $W$ run through all possible words by the split generation result.  On the analytic side, the image $\cF(\cA)$ is split generated by $\cE^V$, which by Orlov's result is equivalent to $D^b(A(\bar\alpha))$ when $W$ runs through all possible words, hence proving the homological mirror symmetry for the split tori.

\brmk We should remark on some technical points regarding the definition of $\cF$ in Abouzaid-Smith's work.  The construction combines several lines of considerations, but one of the key part is that one needs to use a chain-level ($A_\infty$) model $\cF_\infty$, then regard $\cF$ as its reduction to the homotopy level.  Lemma \ref{l:isomorphic} is only a rephrase that $\cF_\infty$ is a quasi-isomorphism to the image.  For complex dimension $1$, the functor on the homological level that we considered essentially follows from \cite{PZ98}.

To define $\cF_\infty$ for higher dimensions, we first recall the $dg$-enhancement of the derived category of $A(\bar\alpha)$, denoted as $D^\pi_\infty(A(\bar\alpha))$.  This is a $dg$-category with essentially the same objects, but the bounded complexes of coherent sheaves should be replaced by their quasi-isomorphic injective resolutions.  Then the morphism groups are morphisms of chain complexes, which are graded by the degree shifts.  A standard procedure shows $D^b_\infty(A(\bar\alpha))$ can be regarded as an $A_\infty$ category obtained by equipping $D^b(A(\bar\alpha))$ with an $A_\infty$-structure with higher operations (i.e. $\mu_i$ with $i\ge3$).  The construction takes a $dg$-enhancement of derived category, then apply the homological perturbation lemma to obtain the minimality (i.e. vanishing $\mu^1$-terms, see \cite{FOOO_Book}\cite{Seidelbook}).

The construction of $\cF_\infty$ then uses the formalism of quilts by considering $L^{W,t}$ as functors from $\fuk(T(\alpha_1,\cdots,\alpha_{n-1}))$ to $\fuk(T(\alpha_n))$.  By induction of homological mirror symmetry from one complex dimension lower, this is quasi-equivalent to $\text{Fun}(D^b_\infty A(\alpha_1,\cdots,\alpha_{n-1}), D^b_\infty A(\alpha_n))$, while the latter is in turn quasi-equivalent to $D^b_\infty A(\bar\alpha)$ from \cite[Theorem 8.15]{To07}.  Strictly speaking, our result looks slightly different from what was proved in \cite{AS10}, but the essential features of the proof already appeared there and were actually used in their proof of Corollary 1.4.




\ermk


\subsection{Homological mirror symmetry of special isogenous tori}\label{s:MSSI}

We finish the proof of Theorem \ref{t:mirror} in this section.  Throughout we focus on the special isogenous tori with one factor, the general case then follows from the product-functor correspondence described in the previous section due to Abouzaid-Smith.  Alternatively, implementing our proof directly to the case of multiple factors is also pretty straightforward.


Consider the special automorphism of $T(\bar\alpha)$ in Defintion \ref{d:symptori}.
This induces a naive $\Z/l$-action on the Fukaya category on $T(\bar\alpha)$ equipped with $\Z$-gradings.  Let the set $\cB_l$ formed by all Lagrangians of the form $L^{W,t}$ specified by the following condition: $t_i=0$ when $w_i=\bl$, and $t_i=r\alpha/l$ when $w_i=\bm$ for some integer $r$ with $0\le r\le l$.  Let $\cA_l$ be the subcategory of twisted complexes and their direct summands in $D^\pi\fuk(T(\bar\alpha)_l)$ formed by $\cB_l$ (we seem to be using a redundant expression of ``subcategory split generated by $\cB_l$", but we want to distinguish these twisted complexes in the derived categories from objects that are only isomorphic to them).  Then $\mathbb{Z}/l\cdot\cB_l$
 forms a collection of $\mathbb{Z}/l$-Lagrangian submanifolds endowed with equivariant brane structures.  From theorem \ref{t:gen}, $\mathbb{Z}/l\cdot\cB$ split generates the equivariant Fukaya category $\fuk(T(\bar\alpha))^{\Z/l}$, hence also the $\Z/l$-invariant part of $\cA_l$. Moreover, it is easy to see that $\cB_l$ verifies the assumption in Corollary \ref{c:downstairsGen}, thus $(\cA_l)^{\Z/l}$ is equivalent to $D^\pi\fuk(T(\bar\alpha)_l)$.


On the analytic side, let $\cB_l^\vee$ be the set of coherent sheaves on $A(\bar\alpha)$ which are mirror of $\cB_l$ under \eqref{e:mirror}, or equivalently, \eqref{e:Afunctor}. Let $\cA_l^\vee:=\cF(\cA_l)$ be the subcategory of $D^b(A(\bar\alpha))$ generated by the set $\cB^\vee_l$.  The fact that $\cA^\vee_l$ being isomorphic to $\cA_l$ (Lemma \ref{l:isomorphic}) implies that the induced $\Z/l$ action is indeed strict rather than coherent, and its invariant part $(\cA^\vee_l)^{\Z/l}$ is split generated by $G\cdot\cB^\vee_l$.  Hence we have the equivalence

\beq\label{e:sympside}(\cA^\vee_l)^{\Z/l}\cong (\cA_l)^{\Z/l}\cong D^\pi\fuk(T(\bar\alpha)_l).\eeq

It is not difficult to relate $(\cA^\vee_l)^{\Z/l}$ to the more commonly used notion of equivariant derived categories.  Consider the free $\Z/l$-action $g^{an}$ on the analytic torus $A(\bar\alpha)$ as in Example \ref{ex:si}.  One checks from \eqref{e:mirror} that the action $g$ and $g^{an}$ are equivariant with respect to the functor $\cF$ on $\cB_l$ and $\cB_l^\vee$ (this partially justifies our claim that $g$ and $g^{an}$ are mirror actions.  One may further check the full generality of the claim by incoporating local systems as in \cite{AS10}\cite{Ha13} and identify points on the analytic tori and the corresponding skyscraper sheaves, but we will not use this here). Therefore, one may talk about the derived category of equivariant sheaves  $D^b_{\Z/l}Coh(A(\bar\alpha))$ (which is equivalent to the equivariant derived category when the group is finite, see \cite{BL94book}).  We regard $D^b_{\Z/l}Coh(A(\bar\alpha))$ as a (non-full) triangulated subcategory of $D^bCoh(A(\bar\alpha))$.  By \cite[Theorem 2.1]{El09} $D^b_{\Z/l}Coh(A(\bar\alpha))$ is again a split closure of $G\cdot\cB^\vee_l$ (regarded as equivariant sheaves) hence is equivalent to $(\cA^\vee)^{\Z/l}$.  The freeness of the $\Z/l$ further identifies the equivariant derived category to the derived category of the quotient (see for example \cite[Example 1.38]{BBHbook}), i.e. $(\cA^\vee)^{\Z/l}\cong D^b_{\Z/l}(A(\bar\alpha))\cong D^b(A(\bar\alpha)_l)$.  Combining this with \eqref{e:sympside} we conclude Theorem \ref{t:mirror}.

\subsection{Derived equivalences between lattice quotients over $\Lambda$}\label{s:dual}

Theorem \ref{t:mirror} reduces our study of Theorem \ref{t:distinguish} to their mirrors.
To understand when $D^b(A(\bar\alpha)_{\bar l})$ is equivalent to $D^b(A(\bar\alpha')_{\bar l'})$. We need to recall how to dualize abelian varieties and morphisms between them.

Recall that the dual of an abelian variety $A$, denoted as $\wh A$, is the identity component of the Picard scheme, i.e. $\wh A=Pic^0(A)$ (see \cite{Mu70}). Given a homomorphism (as group schemes) of abelian variety $\phi:A\rightarrow B$, a \textit{dual homomorphism} $\wh \phi$ is simply defined by the pull-back. A homomorphism $f:A\times \hat A\rightarrow B\times \hat B$ can be represented as a matrix

$$f=\left(\begin{array}{cc} f_1 & f_2\\
                            f_3 & f_4\end{array}\right).$$

Here $f_1:A\rightarrow B$, $f_2:\hat A\rightarrow B$, $f_3:A\rightarrow \hat B$, $f_4:\hat A\rightarrow \hat B $ are again homomorphisms between abelian varieties. We define

$$\wt f=\left(\begin{array}{cc} \wh f_4 & -\wh f_2\\
                            -\wh f_3 & \wh f_1\end{array}\right).$$

Following \cite{Po96} \cite{Or02}, we say $f$ is an \textit{isometric isomorphism} if $f^{-1}=\wt f$. The following general theorem is our main tool:

\bthm[\cite{Or02}\cite{Po96}]\label{t:dereq} Two abelian varieties $A$ and $B$ over a field $k$ with characteristic $0$ are derived equivalent if and only if there is an isometric isomorphism $f:A\times \wh A\rightarrow B\times \wh B$.
\ethm

Our next task is to explicitly write down the implication of Theorem \ref{t:dereq} on an analytic torus, i.e. a lattice quotient of $(\Lambda^*)^n$. First we need to explain how to dualize a lattice quotient $\sT=\Tan/\Gamma$. Let $\Gamma=\ip{e_i}_{i=1}^n$ for $e_i\in(\Lambda^{*})^n$ an integral basis for the lattice, and $\wh \Tan=Hom(\Gamma,\Lambda^*)$. This is a split (analytic) torus with character group natually identified by evaluation on $\Gamma$

\beq\label{e:dualtorus}\begin{aligned} Hom(\Gamma,&\Lambda^*)&\xrightarrow{\sim}&(\Lambda^*)^n\\
           &\vp&\mapsto&(\vp(e_1),\cdots,\vp(e_n))\end{aligned}\eeq

The dual lattice $\wh\Gamma$ is the restriction of the character group $X(T)$ to $\Gamma$, thus under the above identification, can be represented as

$$\wh\Gamma=\ip{(z_i(e_1),\cdots,z_i(e_n))}_{i=1}^n.$$

\bthm[\cite{BoL91}, Theorem 2.1]\label{t:abdual} If $A=\Tan/\Gamma$ is an abelian variety, then
$$\wh{\Tan/\Gamma}=\wh \Tan/\wh\Gamma.$$

\ethm

The last piece of general framework in rigid geometry we need is a uniformization theorem of a morphism between two lattice quotients.

\bthm[\cite{Ge70}, see also Section 5.3, \cite{PaSurvey}]\label{t:abhom} Let $Hom((\Tan_A,\Gamma_A),(\Tan_B,\Gamma_B))$ be a group of homomorphism of the analytic tori sending $\Gamma_A$ to $\Gamma_B$. Then the natural map $Hom((\Tan_A,\Gamma_A),(\Tan_B,\Gamma_B))\longrightarrow Hom(\Tan_A/\Gamma_A,\Tan_B/\Gamma_B)$ is a bijection.

\ethm

In other words, any homomorphism $f:\Tan_A/\Gamma_A\rightarrow \Tan_B/\Gamma_B$ can be lifted uniquely to $\ol f:\Tan_A\rightarrow\Tan_B$ so that $\ol f(\Gamma_A)\subset\Gamma_B$. We will call $\ol f$ the \textit{uniformization of }$f$.

The above theorems enable us to dualize a morphism $h:A\rightarrow B$ between lattice quotients naturally. Namely, there is a natural pairing

$$\ip{\cdot,\cdot}:\wh T\times\Gamma\rightarrow \Lambda^*$$

by evaluation. In particular, this defines a pairing $\wh\Gamma\times\Gamma\rightarrow\Lambda^*$. Let $\ol h:T_A\rightarrow T_B$ be the uniformization of $h$, then $\wh h:\wh B\rightarrow \wh A$ can be described by its uniformization

\beq\label{e:dualmap}\ip{\ov{\wh h}(\wh w), v}=\ip{\wh w,\ol h(v)}\eeq

for $\wh w\in\wh\Gamma_B,v\in\Gamma_A$. Note that since components of $\ol{\wh f}$ are simply group characters, \eqref{e:dualmap} is sufficient to determine $\ol{\wh h}$ due to the fullness of the lattices.

\bex\label{ex:dual} For $A(\bar\alpha)_{l}$ defined in Example \ref{ex:si}, the dual lattice can be computed explicitly as the evaluation of group characters at lattice points in $\Gamma$. This yields

$$\wh\Gamma_l^{\bar\alpha}=\ip{(q^{\alpha_1/l},1,\cdots,1),(q^{\alpha_i/l},1,\cdots,1,q^{\alpha_i},\cdots,1)}_{i=2}^n.$$

We will denote the valuation matrix

$$\wt Q_{l}^{\bar\alpha}=\left(
 \begin{array}{cccc}
   \alpha_1/l\\
   \alpha_2/l&\alpha_2& \text{\huge0} & \\
   \vdots &&\cdots&\\
   \alpha_n/l &\text{\huge0} & & \alpha_n\end{array}\right)$$


\eex

 \hskip 3mm With these general preparations, we may now explain how to reduce the isometric isomorphism condition in Theorem \ref{t:dereq} to an elementary linear algebra condition. By taking the valuation map, we obtain a matrix representation for a lattice $\Gamma$ by listing $(\sigma(e_1),\cdots,\sigma(e_n))$ as column vectors. Denote this matrix by $M_\Gamma$, then the corresponding matrix for $\wh\Gamma$ is precisely its transpose $M_\Gamma^T$.

Suppose $D^b(T_A/\Gamma_A)\cong D^b(T_B/\Gamma_B)$, by Theorem \ref{t:dereq}, \ref{t:abdual} and \ref{t:abhom}, there is a uniformization map

$$\ol f:T_A\times\wh T_A\rightarrow T_B\times\wh T_B,$$

which induces an isomorphism of lattices $\Gamma_A\times\wh \Gamma_A\xrightarrow{\sim}\Gamma_B\times\wh\Gamma_B$. In the rest of the computation we again consider all maps and vectors of the lattices as column vectors in $\R^{n}$ after taking valuations on components.
Let $\Gamma_{A(B)}=\ip{e_1^{A(B)},\cdots,e_n^{A(B)}}$ and $\wh \Gamma_{A(B)}=\ip{\wh e^1_{A(B)},\cdots,\wh e^n_{A(B)}}$. Then

\beq\label{e:matrixAn} \begin{aligned}&\ol f(e^A_i)=\sum_j F_i^j e^B_j+\sum_j G_{ij}\wh e^j_B,\\
&\ol f(\wh e^i_A)=\sum_j H^{ij}e_j^B+\sum_j I^i_j \wh e^j_B, \hskip 3mm \text{for $1\le i,j\le n$}.
\end{aligned}\eeq

We assembly this into a matrix representation

\beq\label{e:matrix1}\ol f(e_1^A,\cdots,e_n^A,\wh e^1_A,\cdots,\wh e^n_A)^T=\left(\begin{array}{cc} F &G\\
                               H &I \end{array}\right)\cdot\left(\begin{array}{c}e_1^B\\
                               \vdots \\
                               e_n^B\\
                               \\
                               \wh e^1_B\\
                               \vdots\\
                               \wh e^n_B\end{array}\right).\eeq

By the isometric assumption, the inverse of $\ol f$ takes the form

\beq\label{e:matrix2}\ol f^{-1}(e_1^B,\cdots,e_n^B,\wh e^1_B,\cdots,\wh e^n_B)^T=\left(\begin{array}{cc} \wh I &-\wh G\\
                                -\wh H& \wh F \end{array}\right)\cdot\left(\begin{array}{c}e_1^A\\
                               \vdots \\
                               e_n^A\\
                               \\
                               \wh e^1_A\\
                               \vdots\\
                               \wh e^n_A\end{array}\right).\eeq

Note that $F,G,H,I, \wh F,\wh G,\wh H,\wh I$ are all integer matrices, and matrices they formed in \eqref{e:matrix1}\eqref{e:matrix2} are in $GL(2n,\Z)$. Recall the pairing between the dual lattice quotients $\ip{\cdot,\cdot}$ from \eqref{e:dualmap}, let $Q_{A(B)}=(\ip{\wh e^i_{A(B)},e_j^{A(B)}})_{1\le i,j\le n}$ ($j$ is the running index for columns). Then we have:

\beq\label{e:matrix3}\begin{aligned}Q_BF^T=\wh FQ_A,\\
                   IQ_B=Q_A\wh I^T,\\
                   GQ_B=Q_A^T\wh G^T,\\
                   HQ^T_B=Q_A\wh H^T.\end{aligned}\eeq

Combining \eqref{e:matrix1}\eqref{e:matrix2}\eqref{e:matrix3}, a straightforward matrix computation shows:

$$\left(\begin{array}{cc} F &G\\
                               H &I \end{array}\right)\cdot \left(\begin{array}{cc} \text{\Large0} &Q_B^T\\
                               -Q_B &\text{\Large0} \end{array}\right)\cdot \left(\begin{array}{cc} F^T &H^T\\
                               G^T &I^T \end{array}\right)=
                               \left(\begin{array}{cc} \text{\Large0} &Q_A^T\\
                               -Q_A &\text{\Large0} \end{array}\right)$$

The reverse direction works equally well. Therefore, we obtain

\blem\label{l:b-side} $D^b(T_A/\Gamma_A)\cong D^b(T_B/\Gamma_B)$ if and only if $\left(\begin{array}{cc} \text{\Large0} &Q_B^T\\
                               -Q_B &\text{\Large0} \end{array}\right)$ and $\left(\begin{array}{cc} \text{\Large0} &Q_B^T\\
                               -Q_B &\text{\Large0} \end{array}\right)$ are congruent by $GL(2n,\Z)$.

\elem

We now consider the case of an special isogenous torus given in Example \ref{ex:si}. For $A=\wt A(\bar\alpha)_{l}$, $Q_A$ is precisely the matrix $\wt Q_{l}^{\bar\alpha}$ in Example \ref{ex:dual}, after taking $\wh e^i$ as the group character $z_i$. Note also that $(\wt Q_{\bar l}^{\bar\alpha})^T=(M^{\bar\alpha}_l)^{-1}$ as defined in Section \ref{s:sympSI}. Therefore, given $\alpha'$, $l'$,

$$\left(\begin{array}{cc} \text{\Large0} &(\wt Q_{l}^{\bar\alpha})^T\\
                               -\wt Q_{l}^{\bar\alpha} &\text{\Large0} \end{array}\right)\sim \left(\begin{array}{cc} \text{\Large0} &(\wt Q_{l'}^{\bar\alpha'})^T\\
                               -\wt Q_{l'}^{\bar\alpha'} &\text{\Large0} \end{array}\right)\Longleftrightarrow
                               \Omega_{l}^{\bar\alpha}\sim \Omega_{l'}^{\bar\alpha'}.$$

Here $\sim$ is the equivalence relation given by congruence class given by $GL(2n,\Z)$. We then conclude the main theorem of this section from Lemma \ref{l:LSeqCong}:

\bthm\label{t:LSeqS} $T(\bar\alpha)_{l}$ is linearly symplectomorphic to $T(\bar\alpha')_{l'}$ if and only if
$A(\bar\alpha)_{l}$ is derived equivalent to $A(\bar\alpha')_{l'}$.

\ethm

Again, since the lattice computations involved carries over to products with no extra complications, Theorem \ref{t:LSeqS} works for any special isogenous symplectic/analytic tori.  Combining \ref{t:mirror} and \ref{t:LSeqS}, one concludes the proof of \ref{t:distinguish}.

\bibliographystyle{plain}

\bibliography{FukRef}

\begin{thebibliography}{10}

\bibitem{AbGen}
Mohammed Abouzaid.
\newblock A geometric criterion for generating the {F}ukaya category.
\newblock {\em Publ. Math. Inst. Hautes \'Etudes Sci.}, (112):191--240, 2010.

\bibitem{AS10}
Mohammed Abouzaid and Ivan Smith.
\newblock {Homological mirror symmetry for the 4-torus}.
\newblock {\em Duke Math. J.}, 2:445--543, 2010.

\bibitem{AuNotes}
Denis Auroux.
\newblock Notes from a topics course on mirror symmetry.
\newblock {\em http://math.berkeley.edu/\~{}auroux/277F09/index.html}, 2009.

\bibitem{BS01}
Paul Balmer and Marco Schlichting.
\newblock Idempotent completion of triangulated categories.
\newblock {\em J. Algebra}, 236(2):819--834, 2001.

\bibitem{BBHbook}
Claudio Bartocci, Ugo Bruzzo, and Daniel Hern{\'a}ndez~Ruip{\'e}rez.
\newblock {\em Fourier-{M}ukai and {N}ahm transforms in geometry and
  mathematical physics}, volume 276 of {\em Progress in Mathematics}.
\newblock Birkh\"auser Boston, Inc., Boston, MA, 2009.

\bibitem{BL94book}
Joseph Bernstein and Valery Lunts.
\newblock {\em Equivariant sheaves and functors}, volume 1578 of {\em Lecture
  Notes in Mathematics}.
\newblock Springer-Verlag, Berlin, 1994.

\bibitem{BC1}
Paul Biran and Octav Cornea.
\newblock Lagrangian cobordism. {I}.
\newblock {\em J. Amer. Math. Soc.}, 26(2):295--340, 2013.

\bibitem{BC2}
Paul Biran and Octav Cornea.
\newblock Lagrangian cobordism and {F}ukaya categories.
\newblock {\em Geom. Funct. Anal.}, 24(6):1731--1830, 2014.

\bibitem{BL04}
Christina Birkenhake and Herbert Lange.
\newblock {\em Complex abelian varieties}, volume 302 of {\em Grundlehren der
  Mathematischen Wissenschaften [Fundamental Principles of Mathematical
  Sciences]}.
\newblock Springer-Verlag, Berlin, second edition, 2004.

\bibitem{BO01}
Alexei Bondal and Dmitri Orlov.
\newblock Reconstruction of a variety from the derived category and groups of
  autoequivalences.
\newblock {\em Compositio Math.}, 125(3):327--344, 2001.

\bibitem{BoL91}
Siegfried Bosch and Werner L{\"u}tkebohmert.
\newblock Degenerating abelian varieties.
\newblock {\em Topology}, 30(4):653--698, 1991.

\bibitem{CH13}
Cheol-Hyun Cho and Hansol Hong.
\newblock {Finite group actions on Lagrangian Floer theory}.
\newblock {\em arXiv:1307.457}.

\bibitem{CHL13}
Cheol-Hyun Cho, Hansol Hong, and Siu-Cheong Lau.
\newblock {Localized mirror functor for Lagrangian immersions, and homological
  mirror symmetry for $P^1_{a,b,c}$}.
\newblock {\em arXiv:1308.4651}.

\bibitem{El09}
A.~D. Elagin.
\newblock Semi-orthogonal decompositions for derived categories of equivariant
  coherent sheaves.
\newblock {\em Izv. Ross. Akad. Nauk Ser. Mat.}, 73(5):37--66, 2009.

\bibitem{El14}
Alexey Elagin.
\newblock {On equivariant triangulated categories}.
\newblock {\em arXiv:1403.7027v1}.

\bibitem{FHS}
Andreas Floer, Helmut Hofer, and Dietmar Salamon.
\newblock Transversality in elliptic {M}orse theory for the symplectic action.
\newblock {\em Duke Mathematical Journal}, 80(1):251--292, 1995.

\bibitem{FP04}
Jean Fresnel and Marius van~der Put.
\newblock {\em Rigid analytic geometry and its applications}, volume 218 of
  {\em Progress in Mathematics}.
\newblock Birkh\"auser Boston, Inc., Boston, MA, 2004.

\bibitem{FOOO_Book}
Kenji Fukaya, Yong-Geun Oh, H.~Ohta, and Kauro Ono.
\newblock {\em Lagrangian intersection {F}loer theory: anomaly and obstruction,
  Part I \& II}, volume~46 of {\em AMS/IP studies in advanced mathematics}.
\newblock American Mathematical Society, Providence, RI, 2009.

\bibitem{Gan13}
Sheel Ganatra.
\newblock {Symplectic cohomology and duality for the wrapped Fukaya category}.
\newblock {\em arXiv:1304.7312}.

\bibitem{Ge70}
L.~Gerritzen.
\newblock \"{U}ber {E}ndomorphismen nichtarchimedischer holomorpher {T}ori.
\newblock {\em Invent. Math.}, 11:27--36, 1970.

\bibitem{Happ}
Dieter Happel.
\newblock {\em Triangulated categories in the representation theory of
  finite-dimensional algebras}, volume 119 of {\em London Mathematical Society
  Lecture Note Series}.
\newblock Cambridge University Press, Cambridge, 1988.

\bibitem{Ha13}
Luis Haug.
\newblock {The Lagrangian cobordism group of $T^2$}.
\newblock {\em arXiv:1310.8056}.

\bibitem{KS02}
Mikhail Khovanov and Paul Seidel.
\newblock Quivers, {F}loer cohomology, and braid group actions.
\newblock {\em J. Amer. Math. Soc.}, 15(1):203--271, 2002.

\bibitem{Mu81}
Shigeru Mukai.
\newblock Duality between {$D(X)$} and {$D(\hat X)$} with its application to
  {P}icard sheaves.
\newblock {\em Nagoya Math. J.}, 81:153--175, 1981.

\bibitem{Mu70}
David Mumford.
\newblock {\em Abelian varieties}.
\newblock Tata Institute of Fundamental Research Studies in Mathematics, No. 5.
  Published for the Tata Institute of Fundamental Research, Bombay; Oxford
  University Press, London, 1970.

\bibitem{Or02}
D.~O. Orlov.
\newblock Derived categories of coherent sheaves on abelian varieties and
  equivalences between them.
\newblock {\em Izv. Ross. Akad. Nauk Ser. Mat.}, 66(3):131--158, 2002.

\bibitem{Or09}
Dmitri Orlov.
\newblock Remarks on generators and dimensions of triangulated categories.
\newblock {\em Mosc. Math. J.}, 9(1):153--159, back matter, 2009.

\bibitem{PaSurvey}
Mihran Papikian.
\newblock Non-{A}rchimedean uniformization and monodromy pairing.
\newblock In {\em Tropical and non-{A}rchimedean geometry}, volume 605 of {\em
  Contemp. Math.}, pages 123--160. Amer. Math. Soc., Providence, RI, 2013.

\bibitem{Po96}
A.~Polishchuk.
\newblock Symplectic biextensions and a generalization of the {F}ourier-{M}ukai
  transform.
\newblock {\em Math. Res. Lett.}, 3(6):813--828, 1996.

\bibitem{PZ98}
Alexander Polishchuk and Eric Zaslow.
\newblock Categorical mirror symmetry: the elliptic curve.
\newblock {\em Adv. Theor. Math. Phys.}, 2(2):443--470, 1998.

\bibitem{RS13}
Alexander Ritter and Ivan Smith.
\newblock The open-closed string map revisited.
\newblock {\em arXiv:1201.5880v3}, 2013.

\bibitem{SeGraded}
Paul Seidel.
\newblock Graded {L}agrangian submanifolds.
\newblock {\em Bull. Soc. Math. France}, 128(1):103--149, 2000.

\bibitem{Se03}
Paul Seidel.
\newblock Homological mirror symmetry for the quartic surface.
\newblock {\em arXiv:0310414}, 2003.

\bibitem{SeidelBias}
Paul Seidel.
\newblock A biased view of symplectic cohomology.
\newblock In {\em Current developments in mathematics, 2006}, pages 211--253.
  Int. Press, Somerville, MA, 2008.

\bibitem{Seidelbook}
Paul Seidel.
\newblock {\em Fukaya categories and {P}icard-{L}efschetz theory}.
\newblock Zurich Lectures in Advanced Mathematics. European Mathematical
  Society (EMS), Zurich, 2008.

\bibitem{SeiG2}
Paul Seidel.
\newblock Homological mirror symmetry for the genus two curve.
\newblock {\em J. Algebraic Geom.}, 20(4):727--769, 2011.

\bibitem{SeAnEq}
Paul Seidel.
\newblock Lagrangian homology spheres in {$(A_m)$} {M}ilnor fibres via {$\Bbb
  C^*$}-equivariant {$A_\infty$}-modules.
\newblock {\em Geom. Topol.}, 16(4):2343--2389, 2012.

\bibitem{Seidelbook2}
Paul Seidel.
\newblock {\em Categorical Dynamics and Symplectic Topology}.
\newblock http://www-math.mit.edu/\~{}seidel/. 2013.

\bibitem{Sher14}
Nick Sheridan.
\newblock Homological mirror symmetry for calabi-yau hypersurfaces in
  projective space.
\newblock {\em Inventiones Mathematicae}.

\bibitem{Sher11}
Nick Sheridan.
\newblock On the homological mirror symmetry conjecture for pairs of pants.
\newblock {\em J. Differential Geom.}, 89(2):271--367, 2011.

\bibitem{Ta71}
John Tate.
\newblock Rigid analytic spaces.
\newblock {\em Invent. Math.}, 12:257--289, 1971.

\bibitem{To07}
Bertrand To{\"e}n.
\newblock The homotopy theory of {$dg$}-categories and derived {M}orita theory.
\newblock {\em Invent. Math.}, 167(3):615--667, 2007.

\end{thebibliography}

\end{document}